\documentclass[11pt]{amsart}
\usepackage{amsmath}
\usepackage{amsthm}
\usepackage{amsfonts}
\usepackage{amssymb}
\usepackage{eucal}
\usepackage[all]{xy}
\usepackage{amsxtra}  
\usepackage[pdftex, breaklinks, linktocpage=true, bookmarksopen=true,bookmarksopenlevel=0,bookmarksnumbered=true]{hyperref}

\hypersetup{
pdftitle = {Canonical integral structures on the de Rham cohomology of curves},
pdfsubject = {Integral structures on the de Rham cohomology of curves},
pdfkeywords = {Integral structure, de Rham cohomology, curve, model, correspondences, rational singularities, arithmetic surface, Grothendieck duality, Artin conductor, efficient conductor, simultaneous resolution of singularities},
pdfauthor = {\textcopyright\ Bryden Cais},
pdfcreator = {\LaTeX\ with package \flqq hyperref\frqq},
colorlinks = {true}
}

\CompileMatrices

\DeclareFontFamily{OT1}{rsfs}{}
\DeclareFontShape{OT1}{rsfs}{n}{it}{<-> rsfs10}{}
\DeclareMathAlphabet{\mathscr}{OT1}{rsfs}{n}{it}

\DeclareMathOperator{\Effcond}{Effcond}
\DeclareMathOperator{\dRcond}{dRcond}
\DeclareMathOperator{\Art}{Art}
\DeclareMathOperator{\Sw}{Sw}

\DeclareMathOperator{\can}{can}

\DeclareMathOperator{\Frac}{Frac}

\DeclareMathOperator{\pr}{pr}
 
\DeclareMathOperator{\Hom}{Hom}

\DeclareMathOperator{\Ext}{Ext}
\DeclareMathOperator{\Gal}{Gal}
\DeclareMathOperator{\GL}{GL}

\DeclareMathOperator{\Spec}{Spec}

\DeclareMathOperator{\coker}{coker}
\DeclareMathOperator{\et}{\acute{e}t}

\DeclareMathOperator{\dR}{dR}

\DeclareMathOperator{\tr}{tr}
\DeclareMathOperator{\Tr}{Tr}

\DeclareMathOperator{\Pic}{Pic}

\DeclareMathOperator{\tors}{tors}

\DeclareMathOperator{\Char}{char}

\DeclareMathOperator{\GD}{GD}
\DeclareMathOperator{\len}{length}
\DeclareMathOperator{\cone}{cone}
\DeclareMathOperator{\Tot}{Tot}

\newcommand*{\R}{\ensuremath{\mathbf{R}}}              
\newcommand*{\Z}{\ensuremath{\mathbf{Z}}}               
\newcommand*{\Q}{\ensuremath{\mathbf{Q}}}                           

\renewcommand*{\P}{\ensuremath{\mathbf{P}}}

\newcommand*{\C}{\mathscr{C}}
     
\newcommand*{\F}{\mathscr{F}}
\newcommand*{\G}{\mathscr{G}}

\newcommand*{\M}{\mathscr{M}}
  
\renewcommand*{\O}{\mathscr{O}}                    
 
\newcommand*{\X}{\mathscr{X}}

\newcommand*{\scrHom}{\mathscr{H}\mathit{om}}

\renewcommand*{\H}{\ensuremath{\mathbf{H}}}

\renewcommand*{\int}{\ensuremath{\mathrm{int}}}

\renewcommand*{\Im}{\ensuremath{\mathrm{im}}}
\newcommand*{\D}{\displaystyle}

\theoremstyle{plain}
  \newtheorem{theorem}{Theorem}
  \newtheorem{proposition}[theorem]{Proposition}
  \newtheorem{lemma}[theorem]{Lemma}
  \newtheorem{corollary}[theorem]{Corollary}

\theoremstyle{definition}
  \newtheorem{definition}[theorem]{Definition}

\theoremstyle{remark}
  \newtheorem{example}[theorem]{Example}
  \newtheorem{remark}[theorem]{Remark}
  \newtheorem{remarks}[theorem]{Remarks}

\include{header}

\numberwithin{theorem}{section}  
\numberwithin{equation}{section}

\begin{document}
\bibliographystyle{plain}

\title{Canonical integral structures on the de Rham cohomology of curves}

\author{Bryden Cais}
\address{Centre de recherches math\'ematiques, Montr\'eal}
\email{bcais@math.mcgill.ca}
\thanks{This work was  partially supported by the NSF grant DMS-0502170 and by a Rackham Predoctoral Fellowship.}


\begin{abstract}
       For a smooth and proper curve $X_K$ over the fraction field $K$ of a discrete valuation 
	ring $R$, we explain (under very mild hypotheses) 
	how to equip the de Rham cohomology $H^1_{\dR}(X_K/K)$ with a {\em canonical integral structure}:
	i.e. an $R$-lattice which is functorial in finite (generically \'etale) $K$-morphisms of $X_K$ and which is 
	preserved by the cup-product auto-duality on $H^1_{\dR}(X_K/K)$.  Our construction of this lattice uses a certain class of normal proper models
	of $X_K$ and relative dualizing sheaves.  We show that our lattice naturally contains
	the lattice furnished by the (truncated) de Rham complex of a regular proper $R$-model of $X_K$ 
	and that the index for this inclusion of lattices is a numerical invariant of $X_K$ (we call it the {\em de Rham conductor}).
	Using work of Bloch and of Liu-Saito, we prove that the de Rham conductor of $X_K$
	is bounded above by the Artin conductor, and bounded below by the Efficient conductor.  We then study how the position of our canonical lattice inside the
	de Rham cohomology of $X_K$ is affected by finite extension of scalars.
\end{abstract}

\maketitle
\section{Introduction}

Let $R$ be a discrete valuation ring with field of fractions $K$ and residue field $k$, and consider a smooth
proper (not necessarily connected) curve $X_K$ over $K$.
The Hodge to de Rham spectral sequence degenerates \cite{DeligneIllusie}, so there is a short exact
sequence of finite dimensional $K$-vector spaces 
\begin{equation}    
	\xymatrix{
		0\ar[r] & H^0(X_K,\Omega^1_{X_K/K})\ar[r] & H^1_{\dR}(X_K/K)\ar[r] & H^1(X_K,\O_{X_K})\ar[r] & 0
	}.\label{hodgefiltration}
\end{equation}
This exact sequence is auto-dual with respect to the cup-product pairing, and is
both contravariantly and covariantly functorial (via pullback and trace) in finite generically \'etale\footnote{Note that generic \'etaleness is needed in order to
have a trace morphism on functions.}
morphisms of smooth curves, so in particular 
is equipped with an action of the $\Z$-algebra of correspondences on $X_K$.      

By an {\em integral structure} on (\ref{hodgefiltration}), we mean a short exact sequence of
free $R$-modules that recovers (\ref{hodgefiltration}) after extending scalars to $K$.  The main aim of this paper
is to equip (\ref{hodgefiltration}) with a {\em canonical} integral structure, i.e. one that is functorial in
both pullback and trace by finite generically \'etale $K$-morphisms of $X_K$, hence stable by the action of correspondences,
and which is self-dual with respect to the cup-product pairing.

Our main result is Theorem \ref{canonicalintegral}, which states that when $X_K$ has an admissible model over $R$ (see Definition \ref{admissibledef}), 
then there is a canonical integral structure on (\ref{hodgefiltration}).  
We remark that $X_K$ has such a model when, for example,
each connected component has a $K'$-rational point for some tamely ramified extension $K'$ of $K$ (so {\em a fortiori} when $X_K$ has good or semistable reduction),
or more generally when $X_K$ has a regular model which is cohomologically flat (see Remark \ref{admiexist} and Proposition \ref{tamept}).
In particular, our theorem holds for modular curves over $p$-adic fields with {\em arbitrary} $p^n$-level,
which in general can be very far from having semistable reduction; see Remark \ref{KatzMazur}.
Furthermore,  we point out that the formation of our canonical integral structure on (\ref{hodgefiltration}) commutes with base change to the closed fiber; i.e.
the reduction of our integral structure modulo a uniformizer of $R$ can be described in terms of the geometry of the closed fiber of {\em any} admissible
model $X$ of $X_K$. 

Let us briefly explain our motivation for studying canonical integral structures in the de Rham cohomology of curves.  Consider
the tower of modular curves $\{X(\Gamma_r)\}_r$ over $\Q_p$ where $\Gamma_r$ is the congruence subgroup $\Gamma(p^r)\cap \Gamma_1(N)$ 
for some fixed integer $N$ coprime to $p$.  
For each $r$, the $p$-adic \'etale
cohomology groups $H^1_{\et}(X(\Gamma_r)_{\overline{\Q}_p},\Z_p)$ are finite $\Z_p$-modules equipped with a natural action of $G_{\Q_p}:=\Gal(\overline{\Q}_p/\Q_p)$.
Since the ``forget the $p^r$-level structure" map $X(\Gamma_r)\rightarrow X(\Gamma_0)$ is Galois with Galois group $\GL_2(\Z/p^r\Z)$, one obtains a $p$-adic
Banach space representation of $\GL_2(\Q_p)$
$$\widetilde{H}^1_{\Q_p}:= \left(\varprojlim_s \varinjlim_r H^1_{\et}(X(\Gamma_r)_{\overline{\Q}_p},\Z_p)/(p^s)\right)\otimes_{\Z_p} \Q_p$$ 
that is moreover equipped with an action of $G_{\Q_p}$.  In \cite{Emerton}, Emerton shows that $\widetilde{H}^1_{\Q_p}$
is an admissible representation of $\GL_2(\Q_p)$ (in the sense of Schneider-Teitlebaum \cite{SchTei}) which encodes 
many deep properties about $p$-adic families of modular forms and their Galois representations.
Note that for each $r$, the (Galois) representations $H^1_{\et}(X(\Gamma_r),\Q_p)$ are de~Rham, so $p$-adic Hodge theory provides 
a natural isomorphism of filtered $\GL_2(\Z/p^r\Z)$-representations
\begin{equation}
	\left(H^1_{\et}(X(\Gamma_r),\Q_p)\otimes_{\Q_p} B_{\dR}\right)^{G_{\Q_p}} \simeq H^1_{\dR}(X(\Gamma_r)/\Q_p)
	\label{BdR}
\end{equation}
With (\ref{BdR}) in mind, it is natural to ask if one can define a de Rham analogue of $\widetilde{H}^1_{\Q_p}$.  This analogue should be an admissible ($p$-adic Banach) $\GL_2(\Q_p)$-representation equipped with a decreasing filtration by admissible $\GL_2(\Q_p)$ sub-representations.  In order to $p$-adically complete, 
one needs an integral version
of $H^1_{\dR}(X(\Gamma_r)/\Q_p)$ and its Hodge filtration that is ``sufficiently functorial" so as to be equipped with a $\GL_2(\Z/p^r\Z)$-action.    
For large $r$, the Galois representations 
$H^1_{\et}(X(\Gamma_r),\Q_p)$ are very far from being semistable in the sense that they become semistable only after a large, wildly ramified extension of $\Q_p$.  Thus, the issue of an appropriate integral structure on $H^1_{\dR}(X(\Gamma_r)/\Q_p)$ is a subtle one; in particular, integral de Rham cohomology seems to be the ``wrong" lattice (i.e. integral structure) to use, for a variety of reasons.  Our lattices have three favorable features: they are stable under the action of correspondences on the curve,
they are auto-dual with respect to cup-product, and their reductions modulo $p$ can be described in terms of the geometry of the closed fiber of any admissible model of $X$. 
For these reasons, it is our hope that the canonical integral structures studied in this paper will enable us to construct the right de Rham analogue of 
$\widetilde{H}^1_{\Q_p}$, and that some ``completed" analogue of the $p$-adic Hodge theory isomorphism $(\ref{BdR})$ will hold.

There are essentially two ways to construct a canonical integral structure on (\ref{hodgefiltration}).  One method proceeds via the N\'eron model
of the Jacobian of $X_K$ and the theory of {\em canonical extensions of N\'eron models} (see \cite[I,\S5]{MM}).  The other method
employs a careful study of a certain class of models of $X_K$, and works directly with the geometry of arithmetic surfaces.
We take the second approach in this paper, which has the advantage that many aspects
of the lattices we obtain seem simpler to analyze via models of curves than via N\'eron models of Jacobians; in particular, we can prove quite a bit more
via this method.  We refer the reader to \cite{CaisNeron} for the first approach, and note that when $X_K$ has a model over $R$ with generically smooth closed fiber,
the two methods yield lattices which are naturally isomorphic (see \cite{CaisNeron}).

We now explain the main idea of our construction.
When $X_K$ has good reduction over $R$, a simple argument with minimal regular models shows that 
the Hodge filtration of the de~Rham cohomology of a smooth proper model over $R$ yields a canonical integral structure on (\ref{hodgefiltration}). 
When $X_K$ does not have good reduction, the (relative) de Rham cohomology of a proper and flat (even regular!) model $X$ can in 
general be rather pathological.
Roughly, there are two problems: first, the sheaf of relative differentials of a proper flat regular model $X$ is not generally a flat $\O_X$-module,
resulting in various torsion phenomena in de Rham cohomology.  Second, and perhaps more seriously, when $X$ is not $R$-smooth, the
sheaf of relative differentials $\Omega^1_{X/R}$ is {\em not} the relative dualizing sheaf, so one can not expect integral de Rham cohomology
to have a good duality theory (i.e. trace maps attached to finite morphisms of curves and cup-product auto-duality). 
To remedy the situation, for a proper flat and normal model $X$ of $X_K$ over $R$,
we will study a certain sub-complex $d:\O_X\rightarrow\omega_{X/R} $ of the de Rham complex of $X_K$, where $\omega_{X/R}$
is the relative dualizing sheaf of $X$ over $R$.  We will show that when $X$ is cohomologically flat over $R$, the hypercohomology of this complex
provides an integral structure on (\ref{hodgefiltration}).  In general, this integral structure may depend on the choice of model $X$, but we will show
that when $X$ is {\em admissible} (Definition \ref{admissibledef}), the resulting integral structure is independent of the choice of admissible model.    
Given such independence, proving that the resulting integral structure is functorial in finite $K$-morphisms of $X_K$ then comes down to the
problem of showing that a finite morphism of curves over $K$ (that have admissible models over $R$) 
may be extended to a {\em finite} morphism of admissible models.  This is established in Theorem \ref{RayGrus} using the flattening techniques of Gruson-Raynaud
\cite{RayGrus}.  We note that the corresponding extension property of finite morphisms is {\em false} if one insists on working with regular models \cite{LiuModel}, \cite{Abhyankar}; thus, it is {\em essential} that we do not require our admissible models to be regular (only normal).

The major technical tool in what we do is Grothendieck duality (\cite{RD}, \cite{GDBC}).  Thus, after a brief study of the geometry and properties
of the class of admissible curves and models in \S\ref{admcurve}--\S\ref{admmodel}, for the convenience of the reader and for lack of an adequate reference
in the literature, we work out in detail Grothendieck duality for relative curves in \S\ref{GDtheory}. 
The heart of this work is \S\ref{GDhypercoho}, where we study the two-term complex $\O_X\rightarrow \omega_{X/R}$ and its hypercohomology.
When $R$ is the localization at some finite place of the ring of integers in a number field and 
$X$ is regular, this hypercohomology occurs in Bloch's work \cite[\S2 ]{Bloch} on the 
Artin conductor of $X/R$.  In \S\ref{Bloch}, assuming that $R$ has perfect residue field and that $X_K$ is geometrically connected,
we use Bloch's results (and subsequent work of Liu-Saito \cite{LiuSaito}) to define a numerical invariant,  the {\em de~Rham conductor} of $X_K/K$
(Definition \ref{dRconddef}), which measures how far our lattice in $H^1_{\dR}(X_K/K)$ is from the lattice furnished by the image of integral de Rham
cohomology, and we prove that this invariant is bounded above by the Artin conductor of $X_K/K$ (see Theorem \ref{dRcondbound}).
We conclude in \S\ref{BaseChng} by addressing the following natural question: if $K'/K$ is any finite separable extension, then what is the relationship
between the canonical integral structure on $H^1_{\dR}(X_{K'}/K')$ and the lattice furnished (via extension of scalars) by the canonical integral structure on
$H^1_{\dR}(X_K/K)$?  When $K'$ is an unramified extension of $K$, these two lattices coincide, while for general (ramified) extensions, they need not
agree (nor is it the case that one lattice necessarily contains the other; see Remarks \ref{latcontainment}).  Theorem \ref{latticeposition} provides
precise estimates for the relative positions of these two lattices inside $H^1_{\dR}(X_{K'}/K')$ in general.

Throughout, our hypotheses and methods are kept as general as possible, so that our results will apply to many curves
of arithmetic interest.  In particular, we do not require $K$ to have characteristic zero.  Moreover, we wish to emphasize
that our constructions work for quite general curves; in particular, we do {\em not} impose any semi-stability hypotheses.
As we have indicated, one of the main applications of our
theory that we have in mind is to the case of the modular curves $X(\Gamma_r)$ over $\Q_p$.  These curves do not have geometrically connected
generic fiber.  Moreover, it is natural in proofs to have to extend the ground field, which may destroy connectivity, and so 
(with the exception of \S \ref{Bloch} where we need geometric connectivity) we do not assume that $X_K$ connected. 

It is a pleasure to thank Spencer Bloch, Matthew Emerton, Dino Lorenzini, Bill Messing, Martin Olsson, Takeshi Saito, and Karen Smith for several helpful 
exchanges and conversations.  I would especially like to thank Brian Conrad, both for advising the author's Ph.D thesis---which was the genesis
of this manuscript---and for many clarifications and corrections on an earlier version.  I heartily thank the referee for 
correcting a number of errors and for suggesting several improvements.

\tableofcontents

\section{Admissible curves}\label{admcurve}
	
We keep the notation from the introduction, and we put $S=\Spec R$.
In this section, we will introduce and study the class of models of $X_K$ over $S$ that we will use to equip the 
de~Rham cohomology of $X_K$ with a canonical integral structure.
	
By a (relative) {\em curve} over $S$ we mean a flat, finite type, and separated $S$-scheme $f:X\rightarrow S$ of pure relative dimension 1 that is normal with 
smooth generic fiber.  For technical convenience, we do not require relative curves to be proper or connected over $S$.  
Note that the non-regular locus of a curve is a finite subset 
of closed points of the closed fiber: when $R$ is excellent this is well-known \cite[$IV_{2}$, 7.8.6 (iii)]{EGA}, and in general---since $X_K$
is smooth---
it follows from the corresponding statement for the base change $X\otimes_R \widehat{R}$ to the excellent completion $\widehat{R}$
of $R$ by \cite[Lemma 2.1.1]{CES} (note that this uses the fact that a local noetherian ring is regular if it admits a regular faithfully flat 
extension \cite[Theorem 23.7]{matsumura})

  In order to define the class of relative curves that we consider, we must first recall 
the definitions of {\em cohomological flatness} and of {\em rational singularities}.
	
Fix a proper relative curve $f:X\rightarrow S$.  
Recall that $f$ is said to be {\em cohomologically flat} (in dimension zero) if the formation of $f_*\O_X$ commutes with arbitrary
base change.  For our purposes, the significance of cohomological flatness is the following:
	
\begin{lemma}\label{cohflat}
	Let $f$ be as above, and let $\F$ be any $S$-flat coherent sheaf on $X$.  Then the formation of $f_*\F$ commutes with 
	arbitrary base change if and only if $H^1(X,\F)$ is torsion-free; in particular, $f$  
	is cohomologically flat in dimension zero if and only if the $R$-module $H^1(X,\O_X)$ is torsion-free.
\end{lemma}

\begin{proof}
	Let $i:X_k\hookrightarrow X$ be the inclusion of the closed fiber and let $\pi\in R$ be a uniformizer.
	Since $\F$ is $S$-flat, we have a short exact sequence
	\begin{equation*}
		\xymatrix{
			0\ar[r] & {H^0(X,\F)\otimes_R k} \ar[r] & {H^0(X_k,i^*\F)} \ar[r]  & {H^1(X,\F)[\pi]}\ar[r] & 0,
		}
	\end{equation*}
	where $H^1(X,\F)[\pi]$ denotes the $\pi$-torsion submodule of $H^1(X,\F)$.  We conclude that
	the formation of $f_*\F$ commutes with the base change $\Spec k\rightarrow S$ (and hence arbitrary 
	base change) if and only if $H^1(X,\F)$ is torsion-free.
\end{proof}
		
By resolution of singularities for excellent surfaces \cite{Lipman2} and descent arguments from $\widehat{R}$ (see \cite[Theorem 2.2.2]{CES})
it follows that $X$ admits a resolution of
singularities; i.e. there exists a proper birational morphism $\rho:X'\rightarrow X$ with $X'$ regular.\footnote{This fact was first proved by Deligne and Mumford
\cite[\S 2]{DM}, using the same descent arguments from $\widehat{R}$ together with results of Abhyankar \cite{AbRes} and of Hironaka.}
The higher direct image sheaf $R^1\rho_*\O_{X'}$ is then supported at the finitely many non-regular points of $X$;
we say that $X$ has {\em rational singularities} if there exists some such resolution with $R^1\rho_*\O_{X'}=0$.  
Trivially, any regular $X$ has rational singularities.
We will make frequent use of the following fact due to Lipman \cite[Proposition 1.2, (2)]{Lipman},
which shows in particular that whether or not $X$ has rational singularities may be checked using {\em any} resolution
of singularities.
		
\begin{proposition}\label{strongrational}
	The scheme $X$ has rational singularities if and only if for every proper birational morphism
	$\rho: X'\rightarrow X$ with $X'$ normal, we have $R^1\rho_*\O_{X'}=0$.
\end{proposition}
	
\begin{remark}\label{derivedrational}
	Since we require $X$ to be normal, an equivalent formulation of rational singularities is that 
	the natural derived category map $\O_{X}\rightarrow \R \rho_* \O_{X'}$ is a quasi-isomorphism
	for every proper birational map $\rho: X'\rightarrow X$ with $X'$ normal.
\end{remark}
			
We can now introduce our primary object of study in this section:
	
\begin{definition}\label{admissibledef}
	A proper relative curve $f:X\rightarrow S$ is {\em admissible} provided $f$ is cohomologically flat and $X$ has rational singularities.
\end{definition}
	
\begin{remark}
	Our definition of admissible generalizes that of Mazur and Ribet \cite[\S2]{MazRib} in that we impose weaker conditions on the possible singularities
	of $X$ and we do not require the closed fiber to be reduced.  Any proper semistable curve over $S$ is admissible, but the class
	of admissible curves is strictly larger than the class of proper semistable curves (see Example \ref{DMstack} and Remark \ref{KatzMazur}).  
\end{remark}
		
We have formulated the notion of admissibility to ensure that every admissible 
curve over $S$ possesses certain properties. For example, in order for the $R$-module $H^1(X, \O_X)$ to 
yield a ``good" integral structure on $H^1(X_K , \O_{X_K})$, it is crucial that it be torsion free
(such freeness will also be essential for later proofs to work).
The hypothesis that $X$ has rational singularities is also critical for our purposes, 
as it enables us to blow up $X$ at closed points of the closed fiber (and then normalize) without changing 
the cohomology $H^1(X, \O_X)$:
				
\begin{proposition}\label{blowupstable}
	Let $f:X\rightarrow S$ be a proper relative curve with rational singularities and let $\rho:X'\rightarrow X$ be any proper birational morphism
	with $X'$ normal.  Then $X'$ has rational singularities, and the canonical pullback map
	${H^1(X,\O_X)}\rightarrow {H^1(X',\O_{X'})}$ is an isomorphism.  In particular, $X'$ is admissible if and only if $X$ is.
\end{proposition}
		
\begin{proof}
	Since $X$ has rational singularities and $X'$ is normal, Proposition 3.2 (ii) of \cite{Artin}
	ensures\footnote{Note that the {\em proof} given there does not require $R$ to be excellent.} that $X'$ has rational singularities.
	Denoting by $f':X'\rightarrow S$ the structural morphism, the exact sequence of terms of low degree
	of the spectral sequence $R^mf_*R^n\rho_*\O_{X'}\implies R^{m+n}f'_*\O_{X'}$ reads
	\begin{equation*}
		\xymatrix{
		0\ar[r] & R^1f_*(\rho_*\O_{X'})\ar[r] & R^1f'_*\O_{X'}\ar[r] & f_* R^1\rho_*\O_{X'}
		}.
	\end{equation*}
	By Proposition \ref{strongrational}, we have $R^1\rho_*\O_{X'}=0$.  
	Thus, since the canonical map $\O_{X}\rightarrow \rho_*\O_{X'}$ is an isomorphism
	(as $\rho$ is proper birational and $X$ is normal), we deduce that the pullback map
	${H^1(X,\O_X)}\rightarrow {H^1(X',\O_{X'})}$ is an isomorphism.
	It now follows from Lemma \ref{cohflat} that $X'$ is admissible if and only if $X$ is.
\end{proof}
		
As both cohomological flatness and rational singularities have somewhat abstract definitions, 
we wish to record some useful criteria that can be used in practice to show that many 
curves of arithmetic interest are admissible (see Example \ref{DMstack} and Remark \ref{KatzMazur}).
				
\begin{proposition}[Raynaud]\label{RaynaudCrit}
	Let $X$ be a connected proper relative $S$-curve.  If $\Char(k)=0$ then $X$ is cohomologically flat.  
		If $\Char(k)=p>0$ then $X$ is cohomologically flat if the greatest common divisor of the 
	geometric multiplicities of the components of the closed fiber is prime to $p$, 
	or {\em a fortiori} if $X$ admits an \'etale quasi-section.
\end{proposition}
		
\begin{proof}
	  By standard arguments (see, for example, the proof of Corollary 9.1.24 in \cite{LiuBook}), we may reduce to the case that $X_K$ is geometrically connected.
	With this assumption we have $f_*\O_{X}\simeq \O_S$, so since $X$ is normal it follows \cite[6.1.6]{RaynaudPic} that
	$X$ satisfies the hypotheses of Raynaud's ``crit\`ere de platitude cohomologique"  \cite[Th\'eor\`eme 7.2.1]{RaynaudPic},
	which immediately implies the proposition.
\end{proof}

\begin{proposition}\label{rationalcirterion}
	Let $f:X'\rightarrow X$ be a finite generically \'etale surjective map of $2$-dimensional excellent schemes that are normal
	and integral, and assume that $X'$ has rational singularities.  If $\rho:Y\rightarrow X$ is a resolution of singularities,
	then $\deg f$ annihilates $R^1\rho_*\O_{Y}$.  In particular, if $\deg f$ is a unit in $\O_X$ then $R^1\rho_*\O_Y=0$
\end{proposition}
	
Note that a resolution of singularities $\rho: Y\rightarrow X$ always exists thanks to Lipman \cite{Lipman2}.
	
\begin{proof}
	Let $g:Y'\rightarrow Y$ be the normalization of $Y$ in the function field of $X'$.  
	As $f$ is generically \'etale, the function field of $X'$ is a separable extension of the function field of $X$, so 
	 the map $g$ is finite.  Setting $n:=\deg g = \deg f$ (generic degree), we 
	thus have a commutative diagram
	\begin{equation*}
		\xymatrix{
			Y'\ar[r]^g\ar[d]_-{\rho'} & Y\ar[d]^-{\rho}\\
			X'\ar[r]_-{f} & X
		}
	\end{equation*}
	of normal 2-dimensional noetherian excellent schemes with $g$ and $f$ finite of generic degree $n$.
	
	Since $Y$ is regular and $Y'$ is normal (hence Cohen-Macaulay, as $Y'$ has dimension 2)
	and the map $g$ is finite, we conclude by
	\cite[Theorem 23.1]{matsumura} that $g$ is flat.  Thus, $g_*\O_{Y'}$ is finite locally free
	over $\O_Y$, and we have a (ring-theoretic) trace map
	$$\tr_g: g_*\O_{Y'}\rightarrow \O_Y$$  
	with the property that the composite map
	\begin{equation}
		\xymatrix{
			{\O_{Y}} \ar[r]^{g^*} & {g_*\O_{Y'}}\ar[r]^-{\tr_g} & {\O_{Y}}
		}\label{compisn}
	\end{equation}
	is multiplication by $n$.  
	
	We have quasi-isomorphisms
	\begin{equation*}
		\R \rho_* \R g_* \O_{Y'} \simeq \R(\rho\circ g)_* \O_{Y'}  = \R(f\circ \rho')_* \O_{Y'} \simeq \R f_* \R \rho'_* \O_{Y'} \simeq f_* \rho'_*\O_{Y'},
	\end{equation*}
	where the last quasi-isomorphism results because $f$ is finite and $X'$ has rational singularities (see Remark \ref{derivedrational}).
	Thus, $\R\rho_* \R g_* \O_{Y'}\simeq \rho_* g_* \O_{Y'}$, so by applying the functor $\R \rho_*$ to the composite map (\ref{compisn})
	and recalling that $g$ is finite (so $g_*\O_{Y'}\simeq \R g_* \O_{Y'}$)
	we deduce that multiplication by $n$ on $\R \rho_*\O_{Y'}$ factors through $\rho_* g_* \O_{Y'}$.
	Applying the functor $H^1$ on complexes gives the result.	
\end{proof}

The following corollary is usually attributed to Elkik \cite[$\mathrm{II}$, Lemme 1]{elkik}:

\begin{corollary}\label{quobyG}
	Let $A$ be an excellent $2$-dimensional regular domain that is essentially of finite type over a noetherian
	domain $R$, and let $G$ be a finite group
	acting on $A$ by $R$-algebra automorphisms.  If $\#G$ is a unit in $A$ then the ring of $G$-invariants
	$A^G$ has rational singularities.
\end{corollary}

\begin{proof}
	It is easy to see that $A^G$ is normal, and our assumptions ensure that $A^G \rightarrow A$ is finite.
	Thus, the corollary follows at once
	from Proposition \ref{rationalcirterion} with $X=\Spec(A^G)$ and $X'= \Spec A$, using Artin's lemma to ensure
	that $\Frac(A)^G=\Frac(A^G)\rightarrow \Frac(A)$ is finite Galois of degree dividing $\#G$.
\end{proof}

\begin{example}\label{DMstack}
	Using Corollary \ref{quobyG}, we can show that many relative $S$-curves of arithmetic interest
	have rational singularities.  Quite generally, let $\X$ be any regular separated Deligne-Mumford stack of finite type and relative dimension
	one over $S$ and suppose that the finite order of the automorphism group at each geometric point of $\X$ is 
	a unit in $R$.  If the underlying coarse moduli
	space $X$ is a scheme, then it is a relative curve over $S$ with all strictly henselian local rings 
	of the form $A^G$ with $A$ regular and $G$ a finite group of $R$-linear automorphisms of $A$
	with order a unit in $R$.  Since whether or not a local ring has a rational singularity can be checked on 
	  a strict henselization (as the formation of a resolution commutes with ind-\'etale base change),  
	  it follows that in this situation $X$ has rational singularities.
\end{example}

\begin{remark}\label{KatzMazur}
	Consider the modular curve $X(N)_{\Q}$ over $\Q$.
	It is smooth and proper, and is the $\Q$-fiber
	of the coarse moduli space $X(N)$ associated to a certain regular Deligne-Mumford stack $\M_{\Gamma(N)}$ that is proper and flat over $\Z$ \cite{ConradKM}.
	  Since the automorphism group of any geometric point must divide $24$   
	 (as the automorphism group of an elliptic curve has order dividing 24 and the automorphism groups at cuspidal geometric points have order dividing 2)
	 we conclude from the example above that $X(N)_{\Z_{(p)}}$ has rational singularities for $p>3$.  
	Furthermore, the cusp $\infty$ is a $\Q$-rational point of $X(N)_{\Q}$, so
	  by Proposition \ref{RaynaudCrit} we see that $X(N)_{\Q}$ has an admissible model
	over $\Z_{(p)}$ for $p>3$.
	
	Similar reasoning shows that the modular curve $X_0(N)_{\Q}$ has an admissible model
	over $\Z_{(p)}$ when $p>3$ (see Proposition 3 of \cite[\S9.7]{BLR}).
\end{remark}

\begin{remark}
	For a proper relative $S$-curve $X$, both the properties of
	  ``being cohomologically flat" and of ``having rational singularities" have
	nice interpretations in terms of the relative Picard functor $\Pic_{X/S}$ of $X$ over $S$.  On the one hand, thanks to 
	Artin \cite[Theorem 7.3]{Artin2} and Raynaud \cite[Proposition 5.2]{RaynaudPic}, the functor $\Pic_{X/S}$ 
	is an algebraic space if and only if $f$ is cohomologically flat.  On the other hand,
	let $\Pic^0_{X/S}$ be the identity component (see \cite[pg. 154]{BLR}) of the relative Picard functor
	and denote by $J^0$ the identity component of the N\'eron model of the Jacobian of $X_K$.  Let us suppose that
	$\Pic^0_{X/S}$ is a scheme
	(which holds, for example, when the greatest common divisor of the geometric multiplicities   
	of the components of the closed fiber of $X$ is 1; see \cite[\S9.4 Theorem 2 (b)]{BLR}).
	Then the canonical map
	$$\xymatrix@1{ {\Pic^0_{X/S}}\ar[r] & {J^0}}$$
	deduced (via the N\'eron mapping property) from the identification $\Pic_{X_K/K}\simeq J_K$ 
	is an isomorphism if and only if $X$ has rational singularities.  We will not need either of these characterizations
	in this paper (but they both play an essential role in the main comparison theorem of \cite{CaisNeron}).
\end{remark}

\section{Admissible models}\label{admmodel}
	
Fix a smooth proper curve $X_K$ over $K$.
By a {\em model} of $X_K$, we mean a proper curve (recall our conventions) $X$ over $S$ whose generic
fiber is isomorphic to $X_K$ (although we will not explicitly keep track of the generic fiber identification).
A model $X$ of $X_K$ is {\em admissible} if the relative curve $X$ is admissible
in the sense of Definition \ref{admissibledef}.
In this section, we study the admissible models of $X_K$.

\begin{lemma}\label{dominated}
Any two admissible models can be dominated by a third.
\end{lemma}  

\begin{proof}
First note that if $C$ is any flat, finite type, and separated $S$-scheme of pure relative dimension 1 having smooth generic fiber then the normalization map
$\widetilde{C}\rightarrow C$ is finite.  Indeed, by smoothness of $C_K$ we see that the base change of the normalization
 $(\widetilde{C})_{\widehat{R}}\rightarrow C_{\widehat{R}}$
is dominated by the normalization of the base change $\widetilde{(C_{\widehat{R}})}\rightarrow C_{\widehat{R}}$, and the latter map is finite as $\widehat{R}$
is excellent.

If $X_1$ and $X_2$ are admissible models of $X_K$,
we let $\Gamma_K$ denote the graph of the generic fiber isomorphism $(X_1)_K\simeq (X_2)_K$ in the fiber product $(X_1)_K\times_K (X_2)_K$. 
Let $\Gamma$ be the closure of $\Gamma_K$ in $X_1\times_S X_2$. 
Then $\Gamma$ has $K$-fiber $\Gamma_K$ because each $X_i$ is $S$-separated.
We take $X_3$ to be the normalization of $\Gamma$.  
The normalization map $X_3\rightarrow \Gamma$ is finite since $\Gamma_K$ is smooth, so the natural maps $X_3\rightrightarrows X_i$ induced by
the projection maps $\pr_i:X_1\times_S X_2\rightrightarrows X_i$ for $i=1,2$ are proper and birational.  Thus, by Proposition \ref{blowupstable},
$X_3$ is an admissible model of $X_K$.  
\end{proof}

\begin{remark}\label{admiexist}
It is natural to ask if every smooth proper curve $X_K$ over $K$ has an admissible model
over $R$.  Quite generally, we claim that $X_K$ has an admissible model over $R$ if and only if {\em every} regular
proper model $X$ of $X_K$ is cohomologically flat (i.e. admissibe).  Indeed, if $X$ is any regular proper model
of $X_K$ and $X'$ is an admissible model of $X_K$, then by resolution of singularities and
the argument of Lemma \ref{dominated}, one knows that there exists a regular model $X''$ of $X_K$ dominating $X$ and 
$X'$.  We conclude from Proposition \ref{blowupstable} that $X$ must be an admissible model of $X_K$.

Thus, if the residue characteristic of $R$ is zero then any regular model of $X_K$ over $R$ is admissible, due
to Proposition \ref{RaynaudCrit}.  On the other hand, if the residue characteristic of $R$ is $p>0$ then by
\cite[9.4.3 a)]{RaynaudPic}, there exist regular proper curves $X$ that are not cohomologically flat.
 Nonetheless, ``most" curves of interest do have an admissible model, as the Proposition below shows.
 \end{remark}
 
\begin{proposition}\label{tamept}
Let $X_K$ be a smooth connected and proper curve over $K$ with a $K'$-rational point for some tamely ramified
extension $K'$ of $K$.   Then every regular proper model of $X_K$ over $R$ is admissible; in particular, $X_K$
has an admissible model over $R$.
\end{proposition}

\begin{proof}
Let $K'$ be {\em any} finite extension of $K$, with ramification index $e$.  If $X_K$ has a $K'$-point, then we claim
that any regular proper model of $X_K$ over $R$ has an irreducible component in the closed fiber of geometric multiplicity dividing $e$.
Proposition \ref{tamept} follows from this, by applying Proposition \ref{RaynaudCrit} (and using the existence of a regular proper model).
To justify the claim, write $v$ for the valuation on $K$ associated to $R$ and choose an extension $v'$ of $v$ to $K'$.  Denote by 
$R'$ the valuation ring of $v'$, so $R\rightarrow R'$ is a local map of discrete valuation rings with ramification index
$e(v'|v) = e$ (note that in general, the integral closure of $R$ in $K'$ need not be local).  
Let $X'$ be a resolution of singularities of the base change $X\times_R R'$, and denote by $\rho: X'\rightarrow X$
the natural map.  As $X'$ is proper, the valuative criterion for properness ensures that the $K'$-point of $X$ extends to an $R'$-point
of $X'$; in other words, the regular $R'$-scheme $X'$ has a section.  This section {\em must} have image in the smooth locus
(by \cite[\S3.1, Proposition 2]{BGR}) so we conclude that there is an irreducible component in the closed fiber of $X'$ with
geometric multiplicity 1.   Letting $\xi'$ be the generic point of this component and writing $\xi=\rho(\xi')$ for the image of $\xi'$ in $X$,
we then have a canonical local map of discrete valuation rings $\O_{X,\xi}\rightarrow \O_{X',\xi'}$.
Fixing an algebraic closure $\overline{k}$ of the residue field of $R$ and applying $\otimes_R \overline{k}$,
we conclude that the length of $\O_{X, \xi}\otimes_R \overline{k}$ divides the length of $\O_{X',\xi'}\otimes_R \overline{k}$.
But by our choice of $\xi$, we know that the length of $\O_{X',\xi'}\otimes_{R'} \overline{k}$ is 1, and it follows that the length
of $\O_{X',\xi'}\otimes_R \overline{k}$ is just the length of $R'\otimes_R k$, which is $e$.  Thus, the closure of $\xi$ in $X$ is an irreducible
component of the closed fiber with geometric multiplicity that divides $e$, giving the claim.
\end{proof}

Now suppose that $f_K:X_K\rightarrow Y_K$ is any finite morphism of smooth proper curves over $K$.  
A standard argument with resolution of singularities and closures of graphs, as in the proof of Lemma \ref{dominated}, shows that
one can always find regular proper models $X$ and $Y$ of $X_K$ and $Y_K$ and a morphism $f:X\rightarrow Y$
recovering $f_K$ on generic fibers.  It is natural to ask if this can be done with $f$ {\em finite} as well.
As Abhyankar showed \cite{Abhyankar}, the answer is negative in general.
The problem has been studied by Liu and Lorenzini \cite[\S 6]{LiuModel}, who show that there are local obstructions in certain cases
to extending $f_K$ to a finite morphism $f$ between regular models, and give a general class of examples of cyclic
\'etale covers $f_K$ that do not extend to finite morphisms of regular models.  
In the positive direction, Liu \cite[Theorem 0.2]{Liu}
has shown that one can always extend $f_K$ to a finite morphism between (possibly non-regular) semistable models after passing to a finite
separable extension $K'$ of $K$.  Using the flattening techniques of Gruson-Raynaud \cite{RayGrus},
when $X_K$ and $Y_K$ have admissible models we can extend $f_K$        
to a finite morphism of admissible models over $R$ (i.e. we avoid passing to an extension of the base field).  More precisely:

\begin{theorem}\label{RayGrus}
	Suppose that $X$ and $Y$ are admissible $S$-curves and $f_K:X_K\rightarrow Y_K$
	is a finite morphism of their generic fibers.  Then there exist admissible models $X'$ of $X_K$
	and $Y'$ of $Y_K$ which dominate $X$ and $Y$ respectively and a finite morphism $f':X'\rightarrow Y'$
	recovering $f_K$ on generic fibers.
\end{theorem}

\begin{proof}
	Let $\Gamma_K\subseteq X_K\times Y_K$ be the graph of $f_K$.  
	As $Y$ is $S$-separated, 
	the closure $\Gamma$ of $\Gamma_K$ in $X\times Y$ has $K$-fiber $\Gamma_K$.  
	Let $\widetilde{\Gamma}$ be the normalization
	of $\Gamma$.  Since the map $\widetilde{\Gamma}\rightarrow \Gamma$
	is finite (see the proof of Lemma \ref{dominated}) it is proper, and clearly $\widetilde{\Gamma}_K\simeq \Gamma_K\simeq X_K$ via the first projection.
	Thus, composing the normalization map with the projection maps, we obtain a proper birational $S$-morphism
	$\widetilde{\Gamma}\rightarrow X$ and an $S$-morphism $f: \widetilde{\Gamma}\rightarrow Y$ that recovers
	$f_K$ on generic fibers.  By Proposition \ref{blowupstable}, $\widetilde{\Gamma}$ is admissible. 

	If $f$ is quasifinite then it is finite and so we are done.  However, this may not be the case.	
	Applying \cite[Th\'eor\`eme 5.2.2]{RayGrus} with $n=1$ to the morphism $f$, which is generically finite,
	we find that there is a blowup $Y''$ of $Y$ such that the strict transform $\widetilde{\Gamma}'$ of $\widetilde{\Gamma}$
	with respect to $Y'' \rightarrow Y$ is quasi-finite over $Y''$. 
	We let $X'$ be the normalization of $\widetilde{\Gamma}'$ and $Y'$ be the normalization of $Y''$; as we saw in the proof
	of Lemma \ref{dominated},
	these normalization maps are finite \cite[$IV_{2}$, 7.8.6 (ii)]{EGA}.
	  The quasi-finite dominant
	map $X'\rightarrow Y''$ necessarily factors through $Y'$ (by the universal property of normalization), so we obtain a quasi-finite $S$-morphism
	$f':X'\rightarrow Y'$.  Now $X'$ and $Y'$ are evidently proper over $S$ and are admissible by Proposition \ref{blowupstable};	
	moreover, the proper quasi-finite morphism $f'$ must be finite, so we are done (see the diagram below).  
	\begin{equation*}
		\xymatrix{
			{X'}\ar[d]\ar@{-->}[drr]^-{f'}                    &                                  &                   \\
			{\widetilde{\Gamma}'}\ar[d]\ar[drr]      &                                   &   {Y'}\ar[d]\\
			{\widetilde{\Gamma}}\ar[d]\ar[drr]^-{f} &                                   & {Y''}\ar[d]\\
			{\Gamma} \ar[r] & {X\times_S Y} \ar[r]_-{\pr_2}\ar[d]_-{\pr_1} & {Y}\ar[d] \\
					     &  {X}\ar[r]                                                        & {S}
		}
	\end{equation*}
\end{proof}

Let $R\rightarrow R'$ be a local extension of discrete valuation rings with fraction field extension $K\rightarrow K'$ and put $S':=\Spec R'$.
Suppose that $X_K$ has an admissible model $X$ over $S$.  By \cite[Lemma 2.1.1]{CES}, if the relative ramification index $e(R'/R)$ is $1$ and the residue 
field extension $k\rightarrow k'$ is separable, then 
$X\times_S {S'}$ is an admissible curve, as such base change preserves regularity and normality (since $X_K$ is smooth), 
and cohomology commutes with flat base change.
However, if $e(R'/R)>1$ or $k'/k$ is not separable
then $X\times_S{S'}$ may no longer be normal, and its normalization need not be admissible.
In general, it is not clear that $X_{K'}$ has an admissible model (but see Remark \ref{admiexist} and Proposition \ref{tamept}), and if it does, how this model is related to $X$.
Nevertheless, we can still say something about the behavior of admissible models with respect to base change:

\begin{proposition}\label{BC}	
	With the notation above, suppose that $X_{K'}$ has an admissible model over $S'$ and that $K\rightarrow K'$ is finite.
	Then there exists an admissible model $Y$ of $X_K$ dominating $X$ such that the normalization of the base change $Y\times_S S'$ is an admissible
	model of $X_{K'}$ over $S'$.  Furthermore, $Y$ may be taken to be regular.
\end{proposition}

\begin{proof}
	Let $X'$ be an admissible model of $X_{K'}$ over $S'$.  As $S'\rightarrow S$ is surjective, there is a canonical isomorphism
	$X'\times_S \Spec K \simeq X_{K'}$.  By Theorem \ref{RayGrus}, there exist admissible models $Z'$ of $X_{K'}$ over $S'$
	and $Z$ of $X_K$ over $S$ dominating $X'$ and $X$, respectively, and a finite morphism $f:Z'\rightarrow Z$ recovering
	the projection $(X')_K\simeq X_{K'}\rightarrow X_K$ over $K$.  Note that the finiteness of $f$ ensures that the normal scheme $Z'$ is
	 the normalization of $Z$ in the function field of $X'$.
	Let $Y\rightarrow Z$ be any proper birational map with $Y$ normal
	(for example, we could take $Y$ to be a resolution of $Z$), and denote
	by $Y'$ the normalization of the base change $Y\times_S {S'}$.  Since $Y'$ is normal and has the same function field as $Z'$, 
	the map $Y'\rightarrow Z$ necessarily factors (uniquely) through $Z'\rightarrow Z$.
	  Due to Proposition \ref{blowupstable}, $Y$ is admissible (over $S$) and $Y'$ is admissible over $S'$.
\end{proof}

\section{The dualizing sheaf and Grothendieck duality}\label{GDtheory}

Recall that a locally finite type map of locally noetherian schemes is said to be {\em Cohen-Macaulay} (which we abbreviate to ``CM") if it is 
flat and all fibers are Cohen-Macaulay schemes.   The composite of two CM maps is again CM, and the base change of a CM map is CM.
Let $X$ and $Z$ be locally noetherian schemes and $f:X\rightarrow Z$ any CM map with pure relative dimension $n$.
Suppose that $Z$ admits a dualizing complex.
By \cite[Theorem 3.5.1]{GDBC}, 
the relative dualizing complex $f^{!}\O_Z$ has a unique nonzero cohomology sheaf, which is in degree $-n$, and we define
\begin{equation}
	\omega_{X/Z}:=H^{-n}(f^{!}\O_Z).\label{dualizingcplxdef}
\end{equation}
We call the sheaf $\omega_{X/Z}$ the {\em $($relative$)$ dualizing sheaf} (for $X$ over $Z$). 
It is flat over $Z$ by \cite[Theorem 3.5.1]{GDBC}, and locally free if and only if the Cohen-Macaulay 
fibers of $f$ are Gorenstein \cite[V, Proposition 9.3, Theorem 9.1]{RD} (which holds, for example, when $X$ is regular).
Thanks to \cite[Theorem 3.6.1]{GDBC}, the formation of $\omega_{X/Z}$ is compatible with
arbitrary base change on $Z$, and the discussion
preceding \cite[Corollary 4.4.5]{GDBC} shows moreover that the formation of $\omega_{X/Z}$ 
is compatible with \'etale localization on $X$.

\begin{remarks}\label{omegaabbrev}
	We will sometimes write $\omega_f$ for $\omega_{X/Z}$, to emphasize its role as {\em relative} dualizing sheaf, or to conserve notation as the case may be.
	In fact, if $Z$ is Gorenstein with finite Krull dimension then $\omega_{X/Z}$ is also a  dualizing complex \cite[V, \S2]{RD} for the abstract scheme $X$.  
	  Indeed, this follows from the fact that $f^{!}$ carries 
	dualizing complexes for $Z$ to dualizing complexes for $X$, which in turn results from Theorem 8.3 and the Remark following it in \cite[V, \S8]{RD} and the local nature of 	being a dualizing complex.  We will need this fact later.
\end{remarks}

We first apply these considerations to the particular situation that $f: X\rightarrow S$ is a relative curve.
Since $S=\Spec R$ with $R$ a discrete valuation ring, $S$ is regular and of finite Krull dimension; it follows from
\cite[V, \S10]{RD} that $\O_S$ is a dualizing complex for $S$.  Moreover, since
$X$ is normal and of dimension 2, it is automatically Cohen-Macaulay
by Serre's criterion for normality.  
In particular, the dualizing complex is a coherent sheaf $\omega_{X/S}$ concentrated in some degree.
By our discussion, the sheaf $\omega_{X/S}$ is flat over $S$ and its formation is compatible
with arbitrary base change on $S$ and \'etale localization on $X$.  

\begin{remark}\label{explicitdualsheaf}
When $X$ is also regular and proper over $S$, there is a more concrete description of $\omega_{X/S}$ that is well-suited to explicit local calculation.
Indeed, one knows that $X$ is projective over
$S$ \cite[Theorem 2.8]{Lichtenbaum}, so we fix a projective $S$-embedding $i:X \rightarrow \P$ with $\P=\P^n_S$.  Since $f:X\rightarrow S$ is a finite type morphism of regular noetherien schemes,
$i$ is a local complete intersection in the sense of \cite[Definition 6.3.17]{LiuBook}.  We may thus consider the {\em canonical sheaf} \cite[Definition 6.4.7]{Liu}, 
defined in terms of the exact sequence
\begin{equation*} 
	\xymatrix{
		0\ar[r] & {\C_{X/\P}}\ar[r] & {i^*\Omega^1_{\P/S}} \ar[r] & {\Omega^1_{X/S}}\ar[r] & 0 
	}
\end{equation*}
as
\begin{equation*}
	\omega_{X/S}^{\can} := \det(\C_{X/\P})^{\vee} \otimes_{\O_{X}} i^*\det(\Omega^1_{\P/S}).
\end{equation*}
By \cite[(2.5.1), Corollary 3.5.2]{GDBC} (or \cite[Theorem 6.4.32]{LiuBook}), there is a canonical isomorphism $\omega_{X/S}\simeq \omega_{X/S}^{\can}$.
Although we do not explicitly need this fact, it will be used implicitly several times in what follows
(mainly to ensure that what we do regarding the dualizing sheaf is compatible with the work of other authors).
\end{remark}

A basic property of the relative dualizing sheaf that will play a crucial role in what we do (and which justifies the name ``dualizing") is the following:
\begin{lemma}\label{sheavesaredual}
	Let $X$ be a relative $S$-curve.  Then the canonical maps
	\begin{equation*}
		\xymatrix{
			{\omega_{X/S}} \ar[r] & {\R\scrHom^{\bullet}_{X}(\O_X[1],\omega_{X/S}[1])}
		}
	\end{equation*}
	and
	\begin{equation*}
		\xymatrix{
			{\O_X[1]} \ar[r] & {\R\scrHom^{\bullet}_{X}(\omega_{X/S},\omega_{X/S}[1])}
		}
	\end{equation*}
	are quasi-isomorphisms.
\end{lemma}

\begin{proof}
	The first isomorphism is obvious (for any $\O_X$-module in the role of $\omega_{X/S}$) since $\scrHom_X(\O_X,\cdot)$ is naturally 
	isomorphic to the identity functor.
	  By viewing
	$\omega_{X/S}$ as $\R\scrHom^{\bullet}_{X}(\O_X[1],\omega_{X/S}[1])$, our assertion for the second map is precisely the claim that the natural double duality map
	\cite[(1.3.20)]{GDBC}
	\begin{equation*}
		\xymatrix{
			{\O_X[1]}\ar[r] & {\R\scrHom^{\bullet}_X(\R\scrHom^{\bullet}_X(\O_X[1],\omega_{X/S}[1]),\omega_{X/S}[1])}
		}
	\end{equation*}
	is a quasi-isomorphism, which is indeed the case as $\omega_{X/S}[1]$ is a dualizing complex for the abstract scheme $X$; see 
	Remarks \ref{omegaabbrev}.
\end{proof}

Deligne and Rapoport \cite[I, \S2.1.1]{DR} call $\omega_{X/S}$ the sheaf of ``regular differentials." 
This terminology is partially justified since $\omega_{X/S}$ is a subsheaf of the sheaf of differential 1-forms
on the generic fiber of $X$.  Indeed, let $i:X_K\hookrightarrow X$ be the inclusion of the smooth 
generic fiber into $X$.  Since $\omega_{X/S}$ is $S$-flat, the natural map $\omega_{X/S}\rightarrow i_* i^*\omega_{X/S}$
is injective.  Moreover, the general theory of the dualizing sheaf provides a canonical isomorphism $i^*\omega_{X/S}\simeq \Omega^1_{X_K/K}$,
so we have an injective map of $\O_{X}$-modules
\begin{equation}
	\omega_{X/S}\hookrightarrow i_*\Omega^1_{X_K/K}.\label{dualizinginclusion}
\end{equation}
We will see later in Proposition \ref{Omegamap} that more is true: in fact, the canonical map $\Omega^1_{X/S}\rightarrow  i_*\Omega^1_{X_K/K}$
uniquely factors through (\ref{dualizinginclusion}).

We now wish to apply Grothendieck duality to proper relative curves and to proper (e.g., finite) morphisms of relative curves.
Quite generally, if $\rho:X\rightarrow Y$ is any proper map of noetherian schemes and $Y$ admits a dualizing complex, 
then Grothendieck's theory provides a natural transformation of functors (see \cite[\S3.4]{GDBC})
\begin{equation}
	\Tr_{\rho}: \R \rho_* \rho^! \rightarrow 1\label{tracenattransform}
\end{equation}
on the bounded below derived category $\mathbf{D}^+_{c}(Y)$.
If $\G^{\bullet}$ is any bounded below complex of abelian sheaves on $Y$ having coherent cohomology,
and $\F^{\bullet}$ any bounded above complex on $X$ with quasi-coherent cohomology, then (\ref{tracenattransform})
yields a natural map
	\begin{equation}
		\xymatrix{
			{\R \rho_* \R\scrHom_X^{\bullet}(\F^{\bullet}, \rho^!\G^{\bullet})}\ar[r] 
				& {\R\scrHom_Y^{\bullet}(\R \rho_* \F^{\bullet},\G^{\bullet})}
		},\label{GD}
	\end{equation}	
which by Grothendieck-Serre duality \cite[Theorem 3.4.4]{GDBC} is a quasi-isomorphism.

In particular, if $Y$ is a regular scheme with finite Krull dimension and $\rho$ has pure relative dimension 1, then for
$\omega_{X/Y}:=H^{-1}(\rho^{!}\O_Y)$ Grothendieck's trace map
(\ref{tracenattransform}) induces a unique trace map
\begin{equation}
	\gamma_f:R^1\rho_*(\omega_{X/Y}) \rightarrow \O_Y\label{topdifftrace}
\end{equation}
whose formation commutes with arbitrary base change $Y'\rightarrow Y$ with $Y'$ another regular scheme with finite Krull dimension \cite[Theorem 3.6.5]{GDBC},
 and we have: 

\begin{theorem}\label{grothendieckduality}
	Let $f:X\rightarrow S$ be any proper relative curve.  Then the natural maps of free $R$-modules
	\begin{equation}
		\xymatrix{
				{H^{0}(X,\omega_{X/S})}\ar[r] &  {H^1(X,\O_X)^{\vee}}
				}
				\quad\text{and}\quad
			\xymatrix{
				{H^{0}(X,\O_X)}\ar[r] &  {H^1(X,\omega_{X/S})^{\vee}}
				}\label{duality0}
	\end{equation}
	induced by cup-product and the trace map $\gamma_f$ are isomorphisms.  Furthermore, there are short exact
	sequences of $R$-modules
	\begin{equation}
		\xymatrix{
			0\ar[r] & {\Ext^1_R(H^1(X,\O_X)_{\tors}, R)} \ar[r] & {H^1(X,\omega_{X/S})} \ar[r] & {H^0(X,\O_{X})}^{\vee} \ar[r] & 0
			}\label{duality1}
	\end{equation}
	and
	\begin{equation}
		\xymatrix{
			0\ar[r] & {\Ext^1_R(H^1(X,\omega_{X/S})_{\tors}, R)} \ar[r] & {H^1(X,\O_X)} \ar[r] & {H^0(X,\omega_{X/S})}^{\vee} \ar[r] & 0
			}.\label{duality2}
	\end{equation}
	
	In particular the natural maps 
	\begin{equation*}
			\xymatrix{
			{H^{1}(X,\omega_{X/S})}\ar[r] & {H^0(X,\O_X)^{\vee}}
			}
			\quad\text{and}\quad
				\xymatrix{
			{H^{1}(X,\O_X)}\ar[r] & {H^0(X,\omega_{X/S})^{\vee}}
			}
	\end{equation*}
	induced by cup-product and the trace map are isomorphisms of free $R$-modules 
	 	if and only if $f$ is cohomologically flat in dimension 0.
\end{theorem}

\begin{proof}
	Setting $\G^{\bullet}=\O_{S}$ and $\F^{\bullet} = \O_X$ (thought of as complexes in degree zero) in (\ref{GD}) and using
	Lemma \ref{sheavesaredual} gives a natural quasi-isomorphism
	\begin{equation*}
		\R\Gamma(X,\omega_{X/S}[1]) \simeq \R \Hom_R^{\bullet}(\R\Gamma(X,\O_{X}), R).
	\end{equation*}
	Applying $H^{-1}$ and using the spectral sequence
	\begin{equation*}
		E_2^{m,n} = \Ext_R^m(H^{-n}(X,\O_X),R) \implies H^{m+n}(\R \Hom^{\bullet}_R(\R\Gamma(X,\O_{X}), R)),
	\end{equation*}
	(whose only nonzero terms occur when $m=0,1$ and $n=0,-1$)
	we deduce that the natural map
	\begin{equation*}
		\xymatrix{
			{H^0(X,\omega_{X/S})}\ar[r] & {\Hom(H^1(X,\O_{X}),R)}
			}
	\end{equation*}
	is an isomorphism.  Similarly, by applying $H^0$ and recalling that the inclusion of the torsion
	submodule $H^1(X,\O_X)_{\tors}\hookrightarrow H^1(X,\O_X)$ induces an isomorphism on $\Ext^1$'s, 
	we obtain (\ref{duality1}).  The same argument with $\F^{\bullet}=\omega_{X/S}$ and $\G^{\bullet}=\O_S$  
	yields (again using the identification of Lemma \ref{sheavesaredual}) the second map in (\ref{duality0}) and the short
	exact sequence (\ref{duality2}).

	To know that the resulting maps 
	coincide with the natural ones obtained via cup product and the trace map $\gamma_f$, one proceeds
	as in the proof of \cite[Theorem 5.1.2]{GDBC}.  Finally, 
	we know that $H^0(X,\O_X)$ and $H^0(X,\omega_{X/S})$ are free 
	since $f$ is flat and $\omega_{X/S}$ is $S$-flat, 
	and since the fibers of $f$ are of pure dimension 1 the $R$-modules $H^i(X,\O_X)$, $H^i(X,\omega_{X/S})$ are both zero for $i>1$ (by the theorem on formal functions). 
	Due to Lemma \ref{cohflat}, the $R$-module $H^1(X,\O_X)$ is free if and only if $f$ is cohomologically flat in dimension zero.
\end{proof}

Now suppose that $\rho:X \rightarrow Y$ is a proper map of relative $S$-curves.  We will work out some consequences
of Grothendieck duality for $\rho$. 
Denoting by $f:X\rightarrow S$ and $g:Y\rightarrow S$ the structural 
morphisms, by definition of the relative dualizing sheaf (\ref{dualizingcplxdef}) we have quasi-isomorphisms
$$\omega_{X/S}[1]\simeq f^!\O_S \quad\text{and}\quad \omega_{Y/S}[1]\simeq g^!\O_S.$$
Since $\rho$ is an $S$-morphism, we have $f=g\circ \rho$, and there is a natural isomorphism
of functors $f^! \simeq \rho^! g^!$ given as in \cite[(3.3.14)]{GDBC}.  Thus, we have a 
quasi-isomorphism
\begin{equation}
	\omega_{X/S}[1]\simeq f^!\O_S \simeq \rho^! g^!\O_S \simeq \rho^!\omega_{Y/S}[1].\label{dualizingrelation}
\end{equation}
and Grothendieck-Serre duality reads
\begin{equation}
	\R \rho_* \R \scrHom^{\bullet}_X(\F^{\bullet},\omega_{X/S}[1]) \simeq \R\scrHom^{\bullet}_Y(\R\rho_*\F^{\bullet},\omega_{Y/S}[1])\label{dualizingidentification}
\end{equation}
for any bounded above complex of sheaves $\F^{\bullet}$ on $X$ with quasi-coherent cohomology.  
By the definition of the duality isomorphism (\ref{GD}),
the map (\ref{dualizingidentification}) is induced by the natural trace map in the derived category
\begin{equation}
	\xymatrix{
		{\R\rho_*\omega_{X/S}}\ar[r] & {\omega_{Y/S}}
		}\label{derivedtrace}
\end{equation}
obtained from (\ref{dualizingrelation}) and Grothendieck's trace map (\ref{tracenattransform}).
If we apply $H^{-1}$ to the composite of (\ref{derivedtrace}) and the natural map $\rho_*\omega_{X/S}\rightarrow \R\rho_*\omega_{X/S}$, 
we get a canonical $\O_Y$-linear {\em trace map} on relative dualizing sheaves
\begin{equation}
	\xymatrix{
		{\rho_*\omega_{X/S}}\ar[r] & {\omega_{Y/S}}
		}.\label{dualizingtracemap}
\end{equation}

We know from Lemma \ref{sheavesaredual} that via the functor $\R\scrHom^{\bullet}_Y(\cdot,\omega_Y[1])$, the complexes $\O_Y[1]$
and $\omega_{Y/S}$ are dual.  Due to (\ref{dualizingidentification}), the same is true of $\R\rho_*\O_X$ and $\R\rho_*\omega_{X/S}$.
This duality furthermore interchanges pullback and trace:

\begin{proposition}\label{pulltracesheaves}
	Let $\rho:X\rightarrow Y$ be a proper map of relative $S$-curves.  
	Then the trace morphism on relative dualizing sheaves $(\ref{dualizingtracemap})$
	is dual $($via the identifications of Lemma $\ref{sheavesaredual}$$)$ to the canonical pullback map 
	$\xymatrix@1{{\O_Y}\ar[r] & {\rho_*\O_X}}$ on functions.  If $\rho$ is in addition finite and generically \'etale, then there exist unique $\O_Y$-linear trace 
	and pullback maps
	\begin{equation}
		\xymatrix{
			{\rho_*: \rho_*\O_X}\ar[r] & {\O_Y}
			}\label{fnstrace}
	\end{equation}		
	\begin{equation}
		\xymatrix{
			{\rho^*:\omega_{Y/S}} \ar[r] & {\rho_*\omega_{X/S}}
		}\label{dualpull}
	\end{equation}
	which are dual and which coincide with the evident trace and pullback maps on the smooth $K$-fibers.  
	In particular, when $X$ and $Y$ are smooth and $\rho$ is finite \'etale, these maps recover the 
	usual pullback and trace maps on functions and differential forms $($cf. Proposition $\ref{Omegamap}$$)$.
\end{proposition}

\begin{proof}
	Applying $\R\scrHom^{\bullet}_Y(\cdot, \omega_Y)$ to the pullback map $\O_Y\rightarrow \rho_*\O_X$ and using (\ref{dualizingidentification})
	gives a map in the derived category
	\begin{equation*}
		\xymatrix{
			{\R\rho_*\R\scrHom^{\bullet}_X(\O_X,\omega_{X/S}) \simeq \R\scrHom^{\bullet}_Y(\R\rho_*\O_X,\omega_{Y/S})} \ar[r] & {\R\scrHom^{\bullet}_Y(\O_Y,\omega_Y)}
		},
	\end{equation*}
	so by Lemma \ref{sheavesaredual} we obtain a map $\R\rho_*\omega_{X/S}\rightarrow \omega_{Y/S}$ which by definition of the map (\ref{GD})
	is none other than the trace map (\ref{derivedtrace}).
	If $\rho$ is generically \'etale, then on $K$-fibers we claim that 
	the trace map (\ref{dualizingtracemap}) is identified with the usual trace map on differential forms.
	  To check this, by working over a Zariski dense open in $Y_K$ over which $\rho_K$ is finite \'etale,
	  we may reduce to the case that $\rho_K$ is finite \'etale, and we wish to show in this context that (\ref{dualizingtracemap})
	coincides with the map
	\begin{equation*}
	 	\xymatrix{
			{{\rho_K}_*\Omega^1_{X_K/K} \simeq {\rho_K}_*\O_{X_K}\otimes \Omega^1_{Y_K/K}} \ar[r] & {\Omega^1_{Y_K/K}}
			}
	 \end{equation*}
	 that sends a section $b\otimes\eta$ to $\Tr(b)\cdot \eta$, where $\Tr:{\rho_K}_*\O_{X_K}\rightarrow \O_{Y_K}$
	 is the usual trace morphism associated to the finite locally free extension of algebras $\O_{Y_K}\rightarrow {\rho_K}_*\O_{X_K}$.
	This follows easily from the construction of Grothendieck's
	trace map in the finite \'etale case; see for instance (2.7.10) and the description of (2.7.37) in \cite{GDBC}, or \cite[III, \S6]{RD}.
	
Suppose now that $\rho:X\rightarrow Y$ is finite and generically \'etale.  We first claim that we have a natural $\O_Y$-linear trace map (\ref{fnstrace}). 
If $\rho$ is in addition {\em flat}, then this is clear, as we take $\rho_*$ to be the usual trace map on functions $\Tr$ as above.
In general, since $Y$ is normal, $\rho$ must be finite flat over the generic points of the special fiber of $Y$, so 
we can find a Zariski open subset $V\subseteq Y$ with scheme-theoretic preimage
$U:=\rho^{-1}V$ such that $\rho_U:U\rightarrow V$ is finite flat and the complements of $V$ in $Y$ and of $U$ in $X$ consist of finitely many points of 
codimension $2$ (necessarily in the closed fibers).  
Since $X$ and $Y$ are normal, we may thus define (\ref{fnstrace}) to be the composite
\begin{equation}
	\xymatrix{
		{\rho_*\O_X} \ar[r]^-{\simeq} & {\rho_*i_*\O_U = j_*{\rho_U}_*\O_U} \ar[r]^-{j_*({\rho_U}_*)} & {j_*\O_V} & \ar[l]_-{\simeq} {\O_Y}
	}.\label{functionstrace}
\end{equation}
Observe that on generic fibers, the trace map (\ref{functionstrace}) recovers the usual
trace morphism associated to the finite flat and generically \'etale morphism $f_K$.  In particular, the preceding construction is independent of the choice of $V$.

To obtain the desired pullback map (\ref{dualpull}) on dualizing sheaves, we apply $\R\scrHom^{\bullet}_Y(\cdot,\omega_{Y/S})$ to the morphism
(\ref{fnstrace}), which yields a derived category map
\begin{equation}
	\xymatrix{
		{\R\scrHom^{\bullet}_Y(\O_Y,\omega_{Y/S})} \ar[r] & {\R\scrHom^{\bullet}_Y(\R\rho_*\O_X,\omega_{Y/S})\simeq \R\rho_*\R\scrHom^{\bullet}_X(\O_X,\omega_{X/S})}
	}.\label{derivedpullback}
\end{equation}
By Lemma \ref{sheavesaredual} we get a derived category map $\omega_{Y/S}\rightarrow \R\rho_*\omega_{X/S}$, and
since $\rho$ is finite (so $\rho_*\rightarrow \R\rho_*$ is a quasi-isomorphism on coherent sheaves) we obtain the
claimed pullback map (\ref{dualpull}) on dualizing sheaves by applying $H^0$.
By construction, the resulting map is dual to the trace map (\ref{fnstrace}) on functions.

Let us see that on $K$-fibers (\ref{dualpull})
really is the usual pullback map on differential forms.  By shrinking $Y$ if need be, we may assume that $\rho_K$ is finite \'etale.
By definition, (\ref{dualpull}) in this case is the map
\begin{equation*}
	\xymatrix{
		{\Omega^1_{Y_K/K}} \ar[r] & {\scrHom_{Y_K}(\O_{Y_K},\Omega^1_{Y_K/K})} \ar[r] & {\scrHom_{Y_K}({\rho_K}_*\O_{X_K},\Omega^1_{Y_K/K}) 
		\simeq {\rho_K}_*\Omega^1_{X_K/K}} 	}
\end{equation*}
that sends a section $\eta$ of $\Omega^1_{Y_K/K}$ to the section $b'\mapsto \Tr(b')\cdot \eta$ of 
${\Hom_{Y_K}({\rho_K}_*\O_{X_K},\Omega^1_{Y_K/K})}$. The Grothendieck duality isomorphism
\begin{equation*}
\xymatrix{
		{{\rho_K}_*\Omega^1_{X_K/K} \simeq {\rho_K}_*\O_{X_K}\otimes \Omega^1_{Y_K/K}} \ar[r] & {\scrHom_{Y_K}({\rho_K}_*\O_{X_K},\Omega^1_{Y_K/K})}
		}
\end{equation*}
in this context is precisely the map sending a section $b\otimes \eta$ of ${\rho_K}_*\O_{X_K}\otimes \Omega^1_{Y_K/K}$
to the section $(b'\mapsto \Tr(bb')\cdot \eta)$ of $\scrHom_{Y_K}({\rho_K}_*\O_{X_K},\Omega^1_{Y_K/K});$ see (2.7.10) and (2.7.9) in \cite{GDBC},
and compare to \cite[III, pg. 187]{RD}.  Thus, the resulting pullback map on dualizing sheaves (\ref{dualpull}) is the map sending 
the section $\eta$ of $\Omega^1_{Y_K/K}$ to the section $1\otimes \eta$ of ${\rho_K}_*\O_{X_K}\otimes \Omega^1_{Y_K/K}\simeq {\rho_K}_*\Omega^1_{X_K/K}$,
and is hence the usual pullback map on differentials as claimed.
\end{proof}

We end this section by recording a useful consequence of rational singularities:   

\begin{proposition}\label{tracechar}
	Let $X$ be a proper relative $S$-curve, and let $\rho:X'\rightarrow X$ be any proper birational morphism
	with $X'$ normal.  If $X$ has rational singularities, then the trace map $\rho_*: \rho_*\omega_{X'/S}\rightarrow \omega_{X/S}$
	of $(\ref{dualizingtracemap})$
	is an isomorphism.  
\end{proposition}

\begin{proof}
	By the characterization of rational singularities in Remark \ref{derivedrational}, the natural map $\O_{X}\rightarrow \R\rho_*\O_{X'}$ is a quasi-isomorphism,
	so substituting $\G^{\bullet} = \omega_{X/S}$ and $\F^{\bullet} = \O_{X'}$ in (\ref{GD}) and using Lemma \ref{sheavesaredual}, we obtain a quasi-isomorphism
	 \begin{equation}
	 	\xymatrix{
		{\R\rho_*\omega_{X'/S}}\ar[r]^{\simeq} &  {\omega_{X/S}}. 
		}\label{theqim}
	 \end{equation}
	 By definition of the duality isomorphism (\ref{GD}), the map (\ref{theqim}) is none other than the trace map (\ref{derivedtrace}).
	Applying $H^0$ to (\ref{theqim}) thus completes the proof.
\end{proof}

\section{Cohomology of regular differentials}\label{GDhypercoho}

Fix a relative curve $f:X\rightarrow S$.   In this section we study the two-term complex of coherent $\O_X$-modules
whose hypercohomology will be a suitable replacement for the relative de Rham cohomology of $X$ over $S$ and,
when $X$ is admissible, will provide the sought-after canonical integral structure on $H^1_{\dR}(X_K/K)$.
In order to define this complex, we need the following well-known fact:

\begin{proposition}\label{Omegamap}
	 Let $i:X_K\hookrightarrow X$ be the inclusion.  Then the canonical pullback map of differentials $\Omega^1_{X/S}\rightarrow i_*\Omega^1_{X_K/K}$
	 factors through the inclusion $\omega_{X/S}\hookrightarrow i_*\Omega^1_{X_K/K}$ of $(\ref{dualizinginclusion})$ via a unique $\O_X$-linear homomorphism 
	\begin{equation*}
		\xymatrix{
			{c_{X/S}: \Omega^1_{X/S}}\ar[r] & {\omega_{X/S}}
			}.
	\end{equation*}
	Moreover, the restriction of $c_{X/S}$ to any $S$-smooth open subset of $X$ is an isomorphism.
\end{proposition}

\begin{proof}
	Uniqueness is clear, as $\omega_{X/S}$ is a subsheaf of $i_*\Omega^1_{X_K/K}$.
	When $X$ is regular, existence follows immediately from Remark
	\ref{explicitdualsheaf} and \cite[Corollary 6.4.13]{LiuBook}.  In general, let $\rho:X'\rightarrow X$ be a resolution of singularities.
	We then define $c_{X/S}$ to be the composite
	\begin{equation*}
	\xymatrix{
		{\Omega^1_{X/S}}\ar[r] & {\rho_*\Omega^1_{X'/S}} \ar[r]^-{\rho_* c_{X'/S}} & {\rho_*\omega_{X'/S}}\ar[r]^-{(\ref{dualizingtracemap})} & \omega_{X/S}
	},
	\end{equation*}
	where the first map is the usual pullback map on differential 1-forms attached to $\rho$.  Since $\rho$ is an isomorphism over the $S$-smooth locus 
	of $X$, we deduce that the same is true for $c_{X/S}$ from the corresponding statement for $c_{X'/S}$.  
\end{proof}

\begin{remark}
	In general, the morphism $c_{X/S}:\Omega^1_{X/S}\rightarrow \omega_{X/S}$ is neither injective nor surjective, and the
	kernel and cokernel of $c_{X/S}$ encode certain arithmetic information about $X$.  Specifically, when the residue field of $R$ is perfect and $X$ is regular, 
	  Bloch shows \cite[Theorem 2.3]{Bloch} 
	that the Euler characteristic of the two-term (generically exact) complex defined by $c_{X/S}$ is equal to $-\Art(X/S)$,
	where $\Art(X/S)$ is the Artin conductor of $X$ over $S$; see \S \ref{Bloch}.
\end{remark}

Thanks to Proposition \ref{Omegamap}, we have a two-term $\O_S$-linear complex (of $\O_S$-flat $\O_X$-modules) concentrated in degrees 0 and 1
\begin{equation}
	\xymatrix{
		{\omega_{X/S}^{\bullet}:= \O_X}\ar[r]^-{d_S} &{\omega_{X/S}}
	},\label{cplxdefn}
\end{equation}
with $d_S$ the composite of the map $c_{X/S}$ and the universal derivation $\O_X\rightarrow \Omega^1_{X/S}$.
Since $c_{X/S}$ is an isomorphism over the $S$-smooth locus of $f$ in $X$, we see that on the generic fiber
of $X$ the complex $\omega_{X/S}^{\bullet}$ is the usual de Rham complex of $X_K$.

\begin{remark}
	Our choice to have $\omega_{X/S}^{\bullet}$ in degrees 0 and 1 reflects our desire to think of it as a replacement for the de Rham complex
	of $X$ over $S$.  This choice will cause a small nuisance when we apply duality theory to $\omega_{X/S}^{\bullet}$, as $f^{!}\O_S=\omega_{X/S}[1]$,
	so we will frequently have to ``shift degrees" in arguments.  We will sometimes do this implicitly, and in any case without further mention; throughout we follow the
	conventions of \cite[\S 1.3]{GDBC} . 
\end{remark}

For use in what follows, let us record some properties of the complex $\omega_{X/S}^{\bullet}$.

\begin{lemma}\label{autodual}
	Let $X$ be a relative $S$-curve.  Then the natural map
	\begin{equation}
	\xymatrix{
		{\omega_{X/S}^{\bullet}} \ar[r] & {\R\scrHom_{X}^{\bullet}(\omega_{X/S}^{\bullet}, \omega_{X/S}[-1])}
		}\label{map1}
	\end{equation}
	 arising from the two maps in Lemma {\rm \ref{sheavesaredual}} is a quasi-isomorphism.  More generally, if $\rho:X\rightarrow Y$
	is any finite morphism of relative curves over $S$, then there is a canonical quasi-isomorphism
	\begin{equation*}
		\R\rho_*\omega_{X/S}^{\bullet} \simeq \R\scrHom^{\bullet}_Y(\R\rho_*\omega_{X/S}^{\bullet},\omega_{Y/S}[-1]).
	\end{equation*}
\end{lemma}

\begin{proof}
		Denote by $f$ and $g$ the structural morphisms $X\rightarrow S$ and $Y\rightarrow S$, respectively, and
		recall our notational convention from Remark \ref{omegaabbrev}.
	  Applying $\R\rho_*$ to (\ref{map1}) and using Grothendieck duality (\ref{dualizingidentification}), we 
	  easily see that the second claim follows from the first.
		
	As for the first claim, the evident exact triangle
	\begin{equation}
		\xymatrix{
			{\omega_{f}[-1]} \ar[r] & {\omega_{f}^{\bullet}} \ar[r] & {\O_X}
			}\label{filbete}
	\end{equation}
	 and the map (\ref{map1}) fit into a morphism of exact triangles
	\begin{equation}
		\xymatrix{
				{\omega_{f}[-1]} \ar[d]\ar[r] & {\omega_{f}^{\bullet}} \ar[d]^-{(\ref{map1})}\ar[r] & {\O_X}\ar[d]\\
			{\R\scrHom^{\bullet}_{X}(\O_X,\omega_{f}[-1])} \ar[r] & {\R\scrHom^{\bullet}_{X}(\omega_{f}^{\bullet},\omega_{f}[-1])} \ar[r] &
			{\R\scrHom^{\bullet}_{X}(\omega_{f}[-1],\omega_{f}[-1])}
				}\label{rhomfilbete}
	\end{equation}
	where the bottom row is obtained from the top by applying $\R\scrHom_X^{\bullet}(\cdot,\omega_{f}[-1])$,
	and the flanking vertical maps are the ones from Lemma \ref{sheavesaredual}.
	These flanking maps are quasi-isomorphisms, by Lemma \ref{sheavesaredual}, so we conclude that the same is true of (\ref{map1}),
	as desired.
\end{proof}

Now we show that the complex $\omega_{X/S}^{\bullet}$ is both covariantly and contravariantly 
functorial in finite generically \'etale morphisms of relative curves:

\begin{proposition}\label{pulltracecomplex}
	Let $\rho:X\rightarrow Y$ be a finite and generically \'etale morphism of relative $S$-curves.  Then there exist pullback and trace maps
	of complexes 
	\begin{equation}
		\xymatrix{
			{\omega_{Y/S}^{\bullet}} \ar[r]^-{\rho^*} &  {\R\rho_*\omega^{\bullet}_{X/S}}
			}\qquad
		\xymatrix{
			{\R\rho_*\omega^{\bullet}_{X/S}} \ar[r]^-{\rho_*} & {\omega^{\bullet}_{Y/S}}
		}\label{pulltracecplx}
	\end{equation}
	that recover the usual pullback and trace maps on de Rham complexes over the $K$-fibers.
	Moreover, these maps are dual with respect to the duality functor $\R\scrHom^{\bullet}_Y(\cdot,\omega_{Y/S}[1])$.
\end{proposition}

\begin{proof}
	Since $\rho$ is finite, the canonical map $\rho_*\omega_{X/S}^{\bullet} \rightarrow \R\rho_*\omega_{X/S}^{\bullet}$
	is a quasi-isomorphism.  Thus, it suffices to construct the desired maps with $\R\rho_*$ replaced by $\rho_*$.
	For pullback, we take the usual pullback on functions in degree $0$ and the morphism (\ref{dualpull}) in degree 1.
	Similarly, for trace we take the trace map on dualizing sheaves (\ref{dualizingtracemap}) in degree 1 and the map (\ref{fnstrace}) in degree 0.	
	That these really yield maps of complexes may be checked at the generic point of the generic fiber of $Y$
	(since the sheaves involved are all flat and $\rho_*$ is left exact, so preserves injections), 
	where---thanks to our assumption that $\rho$ is generically \'etale---it follows from Proposition \ref{pulltracesheaves} and the 
	description on $K$-fibers.
	The duality statement follows easily from Lemma \ref{autodual} and Proposition \ref{pulltracesheaves}.
\end{proof}

We next turn to the study of the hypercohomology of the complex $\omega_{X/S}^{\bullet}$, which when $X$ is admissible will 
provide the desired canonical integral structure on $H^1_{\dR}(X_K/K)$.

\begin{definition}
	Let $f:X\rightarrow S$ be a proper relative curve over $S=\Spec R$.  We define
	$$H^i(X/R):= \H^i(X,\omega^{\bullet}_{X/S}).$$
\end{definition}

Note that $H^i(X/R)$ is nonzero only for $0\le i\le 2$.
Due to Lemma \ref{pulltracecomplex}, for each $i$ the association $X\rightsquigarrow H^i(X/R)$
is both a contravariant functor (via pullback) and a covariant functor (via trace) 
from the category of proper relative $S$-curves with finite generically \'etale morphisms to the category of finite $R$-modules.
Moreover, as the restriction of $\omega_{X/S}^{\bullet}$ to the generic fiber $X_K$ is the usual de Rham complex
of $X_K$, we have natural isomorphisms
\begin{equation*}
	H^i(X/R)\otimes_R K \simeq H^i_{\dR}(X_K/K).
\end{equation*}
When $X$ is cohomologically flat (for example, admissible), the $R$-modules $H^i(X/R)$
are all free and admit a simple description reminiscent of the structure of the de Rham cohomology
of the smooth and proper curve $X_K$:

\begin{proposition}\label{inthodge}
	Suppose that $f:X\rightarrow S$ is a proper and cohomologically flat relative curve.
	Then there are canonical isomorphisms of free $R$-modules of rank one 
	$$H^0(X/R)\simeq H^0(X,\O_X)\quad \text{and}\quad H^2(X/R)\simeq H^1(X,\omega_{X/S})$$
	and a short exact sequence of finite free $R$-modules
	\begin{equation}
		\xymatrix{
			0\ar[r] & {H^0(X,\omega_{X/S})}\ar[r] & {H^1(X/R)} \ar[r] & {H^1(X,\O_X)} \ar[r] & 0
		}\label{integralhodgefil}	
	\end{equation} 
	which is functorial $($both covariant and contravariant$)$ in finite generically \'etale morphisms of proper and cohomologically flat curves,
	and recovers the Hodge filtration of $H^1_{\dR}(X_K/K)$ after tensoring with $K$. 
\end{proposition}

\begin{proof}
Associated to ``la filtration b\^{e}te'' \cite{DeligneHodge2} of $\omega^{\bullet}_{X/S}$ 
is  the usual hypercohomology spectral sequence
$$E_1^{ij} = H^j(X,\omega_{X/S}^i)\implies H^{i+j}(X/R)$$
which we claim degenerates at the $E_1$-stage.  To see this, it suffices to show that the differentials
$$
	\xymatrix@1{
		{H^j(X,\O_X)} \ar[r]^-{d_S} & {H^j(X,\omega_{X/S})}
	}
$$
are zero for $j=0,1$ since $E_1^{ij}=0$ except possibly for $0\le i,j\le 1$.  
Since the Hodge to de Rham spectral sequence for $X_K$ degenerates at the $E_1$-stage \cite{DeligneIllusie},
the image of $d_S$ must be torsion.  However, $H^j(X,\omega_{X/S})$ is torsion free for $j=0$ due to $S$-flatness
of $\omega_{X/S}$, and for $j=1$ by Theorem \ref{cohflat} (since $X$ is cohomologically flat).
It follows that the edge maps 
$$\xymatrix@1{H^0(X/R) \ar[r] &  {H^0(X,\O_X)}}\quad\text{and}\quad \xymatrix{{H^1(X,\omega_{X/S})}\ar[r] & H^2(X/R)}$$
are isomorphisms, and that we have an exact sequence (\ref{integralhodgefil}).   
The flanking $R$-modules in (\ref{integralhodgefil}) are free
due to $S$-flatness of $\omega_{X/S}$ and cohomological flatness of $X$ (Lemma \ref{cohflat}).  Thus, $H^1(X/R)$
is a free $R$-module as well.  By construction and Proposition \ref{pulltracecomplex}, the exact sequence (\ref{integralhodgefil})  
is functorial (both covariantly and contravariantly) in finite generically \'etale $S$-morphisms of proper and cohomologically flat relative curves,
and recovers the Hodge filtration of $H^1_{\dR}(X_K/K)$ after tensoring with $K$.
\end{proof}

To ease notation, we will abbreviate the short exact sequence (\ref{integralhodgefil}) by $H(X/R)$.
Thus, $H(X/R)$ gives an integral structure on the 3-term short exact sequence $H(X_K/K)$
arising from Hodge filtration of $H^1_{\dR}(X_K/K)$ as in (\ref{hodgefiltration}).

Applying the functor $\Hom_R(\cdot,R)$ to $H(X/R)$, we obtain a sequence
$H(X/R)^{\vee}$ of free $R$-modules 
\begin{equation*}
	\xymatrix{
		0\ar[r] & {H^1(X,\O_X)^{\vee}}\ar[r] & {H^1(X/R)^{\vee}} \ar[r] & {H^0(X,\omega_{X/S})^{\vee}}\ar[r] & 0
	}
\end{equation*}
which is also short exact.

\begin{proposition}\label{selfdualityprop}
	Let $f:X\rightarrow S$ be a proper and cohomologically flat relative curve.
	There is a canonical isomorphism of short exact sequences of free $R$-modules
	\begin{equation}
		\xymatrix{
			H(X/R) \ar[r] & H(X/R)^{\vee}
		}.\label{selfduality}
	\end{equation}
\end{proposition}

\begin{proof}
	Applying $\R f_*$ to the quasi-isomorphism of exact triangles (\ref{rhomfilbete}) and using Grothendieck duality for $f$, we get
	a quasi-isomorphism
	\begin{equation}
		\xymatrix{
				 {\R f_*\omega_{f}^{\bullet}} \ar[r]^-{\simeq} & \R f_*\R\scrHom_X(\omega_{f}^{\bullet},\omega_{f}[-1]) \ar[r]\ar[r]^-{\simeq} 
				 & {\R\scrHom^{\bullet}_{S}(\R f_*\omega_{f}^{\bullet},\O_S)[-2]}
				}
	\end{equation}
	that is compatible with the evident filtrations on both sides arising from (\ref{rhomfilbete}).
	The spectral sequence
	$$E_2^{m,n}=\Ext^m_R(\H^{-n}(X,\omega_{f}^{\bullet}),R)\implies H^{m+n}(\R\Hom^{\bullet}_R(\R\Gamma(X,\omega_{f}^{\bullet}),R))$$ 
	is compatible with the filtrations in (\ref{rhomfilbete}), so since 
	$$\Ext^1_R(\H^2(X,\omega_{f}^{\bullet}),R)=0$$
	by Proposition \ref{inthodge}, we get the claimed isomorphism $H^1(X/R)\rightarrow H^1(X/R)^{\vee}$, compatibly
	with the filtrations.
\end{proof}

\begin{remark}
	There is a more ``concrete" description of the isomorphism (\ref{selfduality}).  Indeed, on the one hand, for any two bounded below complexes $\F^{\bullet}$
	and $\G^{\bullet}$ of $\O_X$-modules, abstract nonsense with canonical flasque resolutions provides a bi-functorial map
	\begin{equation}
		\xymatrix{
			{\H^m(X,\F^{\bullet}) \otimes_{\Gamma(X,\O_X)} \H^n(X,\G^{\bullet})} \ar[r] & {\H^{m+n}(\Tot^{\oplus}(\F^{\bullet}\otimes_{\O_X} \G^{\bullet}))}
			}.\label{abstractnonsense}
	\end{equation}
	On the other hand, there is a natural map of complexes
	\begin{equation}
		\xymatrix{
				{\Tot^{\oplus}(\omega^{\bullet}_{X/S} \otimes_{\O_X}\omega^{\bullet}_{X/S})} \ar[r] & {\omega^{\bullet}_{X/S}}
				},\label{wedge}
       \end{equation}
	defined by the canonical maps 
		\begin{equation}
		\xymatrix{
				{\omega_{X/S}^i\otimes_{\O_X} \omega_{X/S}^j} \ar[r] & \omega^{i+j}_{X/S}
				}\label{wedgedef}
       \end{equation}
	which are zero for $i+j > 1$ and are given by multiplication for $i+j\le 1$. 
	Applying (\ref{abstractnonsense}) with $\F^{\bullet}=\G^{\bullet} = \omega_{X/S}^{\bullet}$ 
	and composing with the map on hypercohomology induced by (\ref{wedge}) yields a natural $R$-linear ``cup product" map
	\begin{equation}
		\xymatrix{
			{H^1(X/R) \otimes_R H^1(X/R)}\ar[r] & {H^{2}(X/R)}
			}.\label{finallycup}
	\end{equation}
	That this map recovers the usual cup product after extending scalars to $K$ follows easily from the definition of the maps 
	(\ref{wedgedef}) and the construction of cup products in de Rham cohomology.
	Composing (\ref{finallycup}) with the edge map isomorphism $H^2(X/R)\simeq H^1(X,\omega_{X/S})$ followed by Grothendieck's trace map (\ref{topdifftrace})
	gives an $R$-bilinear pairing on $H^1(X/R)$ which is easily seen to induce a map of short exact sequences $H(X/R)\rightarrow H(X/R)^{\vee}$.
	That this is the {\em same} map as (\ref{selfduality}) may be checked after extending scalars to $K$, where it follows from the classical description
	of Grothendick duality in terms of cup-product and trace (see, for example, the proof of Theorem 5.1.2 in \cite{GDBC}).
\end{remark}

Before getting to the main result of this section, let us address the following natural question.  
Suppose that $R\rightarrow R'$ is a finite local extension of discrete
valuation rings with associated fraction field extension $K\rightarrow K'$ that is separable.  Let $X$ be a proper and cohomologically flat relative curve
over $S':=\Spec R'$.  Considering $X$ as a proper relative curve over $S$, the characterization of cohomological flatness in Lemma \ref{cohflat} shows that
$X$ is also cohomologically flat over $S$.  Note that there is a canonical identification of generic fibers 
$X\otimes_{R'} K' \simeq X \otimes_R K$, so we may speak of {\em the} generic fiber $X_K$ of $X$.  Since $K'/K$ is separable (i.e.~smooth), the canonical surjective map 
$\Omega^1_{X_K/K}\rightarrow \Omega^1_{X_K/K'}$ is an isomorphism and induces a natural isomorphism 
$$H^1_{\dR}(X_K/K)\simeq H^1_{\dR}(X_K/K'),$$
compatibly with Hodge filtrations.  It thus makes sense to ask how the lattices $H^1(X/R)$ and $H^1(X/R')$ are related inside $H^1_{\dR}(X_K/K)$.  

\begin{theorem}
	With the notation above, let $\Delta\subseteq R'$ be the different of the extension $R\rightarrow R'$.
	As short exact sequences of lattices inside the Hodge filtration $H(X_K/K)$, there are canonical $R'$-linear inclusions 
	$\Delta H(X/R)\hookrightarrow H(X/R')$ and $H(X/R')\hookrightarrow H(X/R)$.  Via these inclusions of lattices, $H^0(X,\omega_{X/S'})$ is identified 
	with $\Delta H^0(X,\omega_{X/S})$, and the induced natural maps of $R'$-modules
	\begin{equation}
		\xymatrix@1{ {\frac{\D H^1(X/R')}{\D \Delta H^1(X/R)}} \ar[r]^-{\simeq} & {\frac{\D H^1(X,\O_X)}{\D \Delta H^1(X,\O_X)}}} \quad\text{and}\quad 
	 	\xymatrix@1{ {\frac{\D H^0(X,\omega_{X/S})}{\D \Delta H^0(X,\omega_{X/S})}} \ar[r]^-{\simeq} & {\frac{\D H^1(X/R)}{\D H^1(X/R')}}} 
		\label{natdelmaps}
	 \end{equation}
	 are isomorphisms.
\end{theorem}

\begin{proof}
	Let us denote by $f':X\rightarrow S'$ the structural morphism, and by $g:S'\rightarrow S$ the natural map. 
	By \cite[Theorem 4.3.3]{GDBC}, there is a canonical isomorphism of $\O_{X}$-modules
	  \begin{equation}
	  	 \omega_{X/S} \simeq \omega_{X/S'} \otimes_{\O_X} {f'}^*(g^! \O_S) \label{differentisom}
	  \end{equation}
	Since $R\rightarrow R'$ is finite, we have (canonically)
	\begin{equation*}
	  	 g^! \O_S \simeq \scrHom_S(g_*\O_{S'},\O_S).
	  \end{equation*}
	  Now by definition of the different, the fractional $R'$-ideal $\Delta^{-1}\subseteq K'$ is
	  $$\Delta^{-1}:=\{\alpha\in K'\ :\ \Tr_{K'/K}(\alpha R')\subseteq R\},$$
	  which is mapped isomorphically (as $R'$-module) onto $\Hom_{R}(R',R)$ via the map
	  $\alpha\mapsto \Tr_{K'/K}(\alpha(\cdot))$.  Thus, denoting by $\widetilde{\Delta}$ the ideal sheaf of $\O_{S'}$ associated to the different $\Delta\subseteq R'$,
	 we have a natural isomorphism of $\O_{X}$-modules
	  \begin{equation}
	  	 \omega_{X/S'} \simeq (f')^*(\widetilde{\Delta}) \otimes_{\O_X} \omega_{X/S} = \widetilde{\Delta} \otimes_{\O_{S'}} \omega_{X/S}
		  \simeq \widetilde{\Delta}\cdot \omega_{X/S}\label{differentisom2}
	  \end{equation}
	  where the final isomorphism results from the fact that $\omega_{X/S}$ is $S$-flat (hence $S'$-flat).  
	  The explicit construction of the isomorphism (\ref{differentisom}) given in \cite[Theorem 6.4.9]{LiuBook} shows that
	  on $K$-fibers, the map (\ref{differentisom2}) is the inverse of the canonical isomorphism $\Omega^1_{X_K/K}\xrightarrow{\simeq} \Omega^1_{X_K/K'}$.   
	  Choosing a generator $\delta$ of $\Delta$, since $\omega_{X/S}$ is $S'$-flat we therefore have a commutative diagram of $\O_X$-modules
	  \begin{equation*}
	  	\xymatrix{
			{i_*\Omega^1_{X_K/K}} \ar[r]^-{\can\circ(\cdot\delta)} & {i_*\Omega^1_{X_K/K'}}  & \ar[l]^-{\simeq}_-{\can} {i_*\Omega^1_{X_K/K}} \\
			{\omega_{X/S}} \ar@{^{(}->}[u]\ar[r]_-{\cdot\delta}^-{\simeq} & {\omega_{X/S'}} \ar@{^{(}->}[u]\ar@{^{(}->}[r] & {\omega_{X/S}}\ar@{^{(}->}[u]
		}
	  \end{equation*}
	  where $i:X_K\hookrightarrow X$ is the inclusion of the generic fiber and $\can$ is the canonical isomorphism.  We claim that this diagram extends
	  to a commutative diagram of two-term complexes:
	  \begin{equation}
	  \xymatrix@R=6pt@C=8pt{ 
		 {i_*\O_{X_K}} \ar[rr]\ar[dr] & & {i_*\O_{X_K}} \ar[dr]  \ar@{=}[rr] & & {i_*\O_{X_K}} \ar[dr] 
		\\ 
		 & {i_*\Omega^1_{X_K/K}} \ar[rr]  & & {i_*\Omega^1_{X_K/K'}}   \ar@{-}[ul]  & & \ar[ll] {i_*\Omega^1_{X_K/K}} 
		\\ 
		 {\O_X} \ar'[r][rr]\ar[dr]\ar@{^{(}->}[uu] & & {\O_X} \ar[dr] \ar@{=}'[r][rr] \ar@{^{(}->}'[u][uu]& & {\O_X} \ar[dr] \ar@{^{(}->}'[u][uu]
		\\
		&{\omega_{X/S}} \ar@{^{(}->}[uu] \ar[rr] & & {\omega_{X/S'}} \ar@{^{(}->}[uu] \ar@{^{(}->}[rr] & & {\omega_{X/S}} \ar@{^{(}->}[uu]
	}\label{cube}	
	\end{equation}  
	Indeed, the commutativity of the left cube is clear, as all the horizontal arrows are (essentially) multiplication by $\delta$.  As for the commutativity of the right cube, 
	it suffices to show that the ($\O_{S'}$-linear) derivation 
	\begin{equation*}
		\xymatrix{
			{\O_X} \ar[r]^-{d_{S'}} & {\omega_{X/S'}} \ar@{^{(}->}[r] & {\omega_{X/S}}
		}
	\end{equation*}
	coincides with $d_{S}$.  Thanks to the flatness of $\omega_{X/S}$ and $\omega_{X/S'}$, this may be checked on generic fibers, where it
	follows from the fact that the canonical isomorphism $\Omega^1_{X_K/K}\rightarrow \Omega^1_{X_K/K'}$ carries $d_K$ to $d_{K'}$.
	Taking hypercohomology of (\ref{cube}), we deduce the following commutative diagram of short exact sequences of lattices in $H(X_K/K)$:
	\begin{equation} 
		\xymatrix{
			0\ar[r] & {H^0(X,\omega_{X/S})} \ar[r]\ar@{^{(}->}[d]^-{\cdot\delta} & {H^1(X/R)}\ar@{^{(}->}[d]^-{\cdot\delta} \ar[r] & 
			{H^1(X,\O_X)}\ar[r]\ar@{^{(}->}[d]^-{\cdot\delta} & 0\\
			0\ar[r] & {H^0(X,\omega_{X/S'})} \ar[r]\ar@{^{(}->}[d] & {H^1(X/R')} \ar[r]\ar@{^{(}->}[d] & {H^1(X,\O_X)}\ar[r]\ar@{=}[d] & 0\\
			0\ar[r] & {H^0(X,\omega_{X/S})} \ar[r] & {H^1(X/R)} \ar[r] & {H^1(X,\O_X)}\ar[r] & 0\\
		}\label{hypdia}
	\end{equation}
	Since $H^1(X,\omega_{X/S})$ is torsion free, we have
	\begin{equation*}	
		H^0(X,\omega_{X/S'}) \simeq H^0(X, \Im(\cdot\delta: \omega_{X/S}\rightarrow \omega_{X/S}))
		\simeq \Im(\cdot\delta: H^0(X,\omega_{X/S})\rightarrow H^0(X,\omega_{X/S}))
	\end{equation*}
	as lattices in $H^0(X,\Omega^1_{X_K/K})$, so the top left (vertical) arrow in (\ref{hypdia}) is an isomorphism.  It then follows at once from (\ref{hypdia}) that the natural maps
	(\ref{natdelmaps}) are isomorphisms.
\end{proof}

In general, the integral structure $H(X/R)$ on the Hodge filtration exact sequence $H(X_K/K)$
of the generic fiber 
may depend on the choice of proper and cohomologically flat model $X$ of $X_K$.
Remarkably, when $X_K$ has an {\em admissible} model, the integral structure $H(X/R)$ for admissible models $X$ is intrinsic to $X_K$: it does not depend on the choice of admissible model $X$.  It is moreover {\em functorial in $K$-morphisms of $X_K$}, and so truly deserves to be called a canonical integral structure.

\begin{theorem}\label{canonicalintegral}
	  Keep the notation introduced immediately prior to Proposition \ref{selfdualityprop}, and
	suppose that $X_K$ is a proper smooth curve over $K$ with an admissible model over $R$.
	\begin{enumerate}
		\item \label{first} For any two admissible models $X$ and $X'$ of $X_K$, there 
		is a unique isomorphism of short exact sequences of finite free $R$-modules
			$$H(X/R)\simeq H(X'/R)$$
		respecting the $K$-fiber identifications.
		\item \label{second} For any finite and generically \'etale morphism $\rho:X_K\rightarrow Y_K$ of proper smooth 
			curves over $K$ that admit admissible
			models $X$ and $Y$, there are canonical $R$-linear homomorphisms of integral structures 
			\begin{equation*}
				\xymatrix{
					H(Y/R) \ar[r]^-{\rho^*} & H(X/R)
				}\quad\text{and}\quad
				\xymatrix{
				H(X/R) \ar[r]^-{\rho_*} & H(Y/R)
				}
			\end{equation*}
			which recover the natural $K$-linear pullback and trace maps induced by $\rho_K$ on the Hodge filtration exact sequences
				$H(X_K/K)$ and $H(Y_K/K)$  after extending scalars to $K$.  
		\item\label{third}	
			The morphisms $\rho^*$ and $\rho_*$ in   $(\ref{second})$ are
			adjoint with respect to the pairing furnished by $(\ref{selfduality})$.  More precisely, under the identifications
			of $H(X/R)$ and $H(Y/R)$ with their $R$-linear duals as in Proposition $\ref{selfdualityprop}$, the morphisms
			\begin{equation*}
				\xymatrix{
					H(X/R)^{\vee} \ar[r]^-{(\rho^*)^{\vee}} & H(Y/R)^{\vee}
				}
					\quad\text{and}\quad		
				\xymatrix{
				H(X/R) \ar[r]^-{\rho_*} & H(Y/R)
				}
		\end{equation*}
		coincide, as do			
			\begin{equation*}
				\xymatrix{
					H(Y/R)^{\vee} \ar[r]^-{(\rho_*)^{\vee}} & H(X/R)^{\vee}
				}
					\quad\text{and}\quad		
				\xymatrix{
				H(Y/R) \ar[r]^-{\rho^*} & H(X/R)
				}.
		\end{equation*}
	\end{enumerate}
\end{theorem}

\begin{proof}
	(\ref{first}).  Uniqueness is obvious, by $R$-freeness.  Suppose $X$ and $X'$ are admissible models of $X_K$.
	Since any two such models can be dominated by a third thanks to Lemma \ref{dominated}, we may suppose there is a (necessarily unique)
	proper birational map $\rho:X'\rightarrow X$ inducing the identity on $K$-fibers.
	Since $\rho$ is birational and $X$ is normal, the natural map $\O_X\rightarrow \rho_*\O_{X'}$ is an isomorphism.
	Moreover, Proposition \ref{tracechar} shows that the trace map on relative dualizing complexes $\rho_*\omega_{X'/S}\rightarrow \omega_{X/S}$
	is an isomorphism, so we obtain a diagram
	\begin{equation*}
		\xymatrix{
			{\O_X} \ar[r]\ar[d]_-{\simeq} & \omega_{X/S}\\
			\rho_*\O_{X'} \ar[r] & \rho_*\omega_{X'/S}\ar[u]_-{\simeq}
		}
	\end{equation*}
	which commutes because this may be checked over $K$, where it is the canonical isomorphism induced by the identification $X_K\simeq X'_{K}$.
	We thus obtain a map of complexes 
	\begin{equation*}
		\xymatrix{
			\omega_{X/S}^{\bullet}\ar[r]^-{\simeq} & \rho_*\omega_{X'/S}^{\bullet}\ar[r] & \R \rho_*\omega_{X'/S}^{\bullet}
			}
	\end{equation*}
	that respects the canonical filtrations on both sides, and hence (upon applying $\R\Gamma(X,\cdot)$) gives a commutative diagram
	\begin{equation}
		\xymatrix{
			0\ar[r] & H^0(X,\omega_{X/S})\ar[r]\ar[d] & H^1(X/R)\ar[r]\ar[d] & H^1(X,\O_X)\ar[r]\ar[d] & 0\\
				0\ar[r] & H^0(X',\omega_{X'/S})\ar[r] & H^1(X'/R) \ar[r] & H^1(X',\O_{X'})\ar[r] & 0
		}\label{compare}
	\end{equation}
	recovering the identification of the Hodge filtrations of $H^1_{\dR}(X_K/K)$ and $H^1_{\dR}(X'_K/K)$ after extending scalars to $K$.
	
	By Proposition \ref{blowupstable},
	the right vertical map of (\ref{compare}) is an isomorphism.  
	By definition, the left vertical map is the inverse of the map on global sections over $X$ induced from the isomorphism
	$\rho_*\omega_{X'/S}\rightarrow \omega_{X/S}$ of Proposition \ref{tracechar}, and is hence an isomorphism too.  
	It follows that all three vertical maps in (\ref{compare}) are isomorphisms.

	(\ref{second}).  Let $\rho:X_K\rightarrow Y_K$ be a finite and generically \'etale morphism of proper smooth curves over $K$ that have 
	admissible models $X$ and $Y$ over $S$.  By Theorem \ref{RayGrus}, there exist admissible models $X'$ of $X_K$ and $Y'$
	of $Y_K$ that dominate $X$ and $Y$, respectively, and a {\em finite} (necessarily generically \'etale) morphism $\rho':X'\rightarrow Y'$ 
	 recovering $\rho$ on generic fibers.
	By Proposition \ref{inthodge}, there exist $R$-linear morphisms of short exact sequences 
	\begin{equation*}
				\xymatrix{
					H(Y'/R) \ar[r]^-{(\rho')^*} & H(X'/R)
				}\quad\text{and}\quad
								\xymatrix{
					H(Y'/R) &\ar[l]_-{(\rho')_*}  H(X'/R)
				}
	\end{equation*}
	recovering the canonical pullback and trace maps, respectively, after extending scalars to $K$.
	Since $Y'$ dominates $Y$ and $X'$ dominates $X$, part (\ref{first}) provides natural isomorphisms
	of short exact seqences
	\begin{equation*}
				H(Y'/R) \simeq H(Y/R) \quad\text{and}\quad H(X'/R)\simeq H(X/R)
	\end{equation*}
	in a manner recovering the canonical isomorphisms induced by the $K$-fiber identifications $Y'_K\simeq Y_K$
	and $X'_K\simeq X_K$ after extending scalars to $K$. We thus obtain the desired maps of integral structures.
	
	(\ref{third}).  The claimed duality of $\rho_*$ and $\rho^*$ is a straightforward consequence of Proposition \ref{pulltracecomplex}
	and basic compatibilities in Grothendieck duality, as we explain below.
	Denote by $g:Y\rightarrow S$ and $f:X\rightarrow S$ the structural morphisms, and for ease of notation, write
	   $$(\cdot)_X^{\vee},\qquad (\cdot)_Y^{\vee}\quad\text{and}\quad (\cdot)_S^{\vee}$$
		for the functors
	  $$\R\scrHom_X(\cdot,\omega_f[-1]),\qquad \R\scrHom_Y(\cdot,\omega_{g}[-1])\quad\text{and}\quad \R\scrHom_S(\cdot,\O_S)$$
	   respectively.  Moreover, if $\bullet$ stands for any one of $g$, $f$, $\rho$, we will denote by
	$\GD_{\bullet}$ the Grothendieck duality isomorphism (\ref{GD}) attached to $\bullet$.  Letting $\rho_*$ and $\rho^*$
	be the pullback and trace morphisms as in (\ref{pulltracecplx}), we claim that there is a commutative diagram of natural 
	quasi-isomorphisms, compatible with filtrations,
	\begin{equation*}
		\xymatrix@C+20pt{
			{(\R g_*\omega_{g}^{\bullet})_S^{\vee}[-2]} \ar[r]^-{ (\R g_*(\rho_*))^{\vee}_S} & 
			{(\R g_*\R\rho_*\omega_{f}^{\bullet})_S^{\vee}[-2]} \ar[r] & {(\R f_*\omega_{f}^{\bullet})_S^{\vee}[-2]}  \\
			{\R g_*(\omega_{g}^{\bullet})_Y^{\vee}} \ar[r]_-{\R g_*(\rho_*)_Y^{\vee}}\ar[u]^-{\GD_g} & 
			{\R g_* (\R\rho_*\omega_{f}^{\bullet})_Y^{\vee}} \ar[u]_-{\GD_g} & \\ 
			 & {\R g_* \R\rho_*(\omega_{f}^{\bullet})_X^{\vee}} \ar[r]\ar[u]_-{\R g_*(\GD_{\rho})} & 
			 {\R f_*(\omega_{f}^{\bullet})_X^{\vee}}\ar[uu]_-{\GD_f}\\
			{\R g_* \omega_{g}^{\bullet}}\ar[r]_-{\R g_*(\rho^*)}\ar[uu]^-{\R g_*(\ref{map1})} & {\R g_* \R\rho_* \omega_{f}^{\bullet}} \ar[r]\ar[u]_-{\R g_*\R\rho_*(\ref{map1})} 
			& {\R f_*\omega_{f}^{\bullet}}\ar[u]_-{\R f_*(\ref{map1})}
		}
	\end{equation*}
	Here, the horizontal maps in the right rectangle are induced by the canonical isomorphism $\R g_*\R\rho_*\simeq \R f_*$
	attached to the composite $f=g\circ\rho$.  The commutativity of the lower left rectangle follows from Proposition \ref{pulltracecomplex},
	by applying $\R g_*$.  The commutativity of the upper left rectangle is due to the functoriality of Grothendieck's duality isomorphism $\GD_g$.
	The lower right rectangle commutes because the isomorphism of functors $\R g_*\R\rho_*\simeq \R f_*$ is natural.  Finally, the upper right
	rectangle commutes because Grothendieck duality is compatible with composition (see \cite[\S3.4]{GDBC}, especially Lemma 3.4.3).
	
	Applying $H^1$ to the large rectangle formed by the outside edges (and using the proof of Proposition \ref{selfdualityprop}) we deduce the commutativity 
	of the diagram
	\begin{equation*}
		\xymatrix{
			{H^1(Y/R)^{\vee}} \ar[r]^-{(\rho_*)^{\vee}} & {H^1(X/R)^{\vee}}\\
			{H^1(Y/R)} \ar[u]_-{\simeq}^-{(\ref{selfduality})} \ar[r]_-{\rho^*} & {H^1(X/R)}\ar[u]^-{\simeq}_-{(\ref{selfduality})}
		}
	\end{equation*}
	compatibly with filtrations, as claimed.
\end{proof}

\section{Comparison with integral de Rham cohomology}\label{Bloch}

Fix a proper relative curve $f:X\rightarrow S$.  
Besides the complex $\omega_{X/S}^{\bullet}$, there 
is another complex of differentials that provides a lattice in the $K$-vector 
space $H^1_{\dR}(X_K/K)$: the 2-term ``truncated de~Rham complex" (in degrees 0,1)
\begin{equation}
	\Omega_{X/S,\le 1}^{\bullet}:=\xymatrix@1{{\O_X}\ar[r]^-{d} & {\Omega^1_{X/S} }}.\label{omegatruncated}
\end{equation}
When $X$ is $S$-smooth, (\ref{omegatruncated}) and $\omega_{X/S}^{\bullet}$ are canonically ismorphic.
In general, thanks to Proposition \ref{Omegamap}, one has a map of complexes
\begin{equation}
	\xymatrix{
		 {\Omega^{\bullet}_{X/S,\le 1}} \ar[r]^{\rho} & {\omega^{\bullet}_{X/S}}
	}\label{rhomap},
\end{equation}
and the discrepancy between the (hyper) cohomologies
of these complexes encodes important arithmetic information about $X$ and $X_K$.

\begin{remark}
	The reader may wonder why we do not also consider the relative de Rham cohomology of $X$ over $S$, given by the
	hypercohomology of the 3-term de Rham complex 
	$$\Omega_{X/S}^{\bullet}:=\xymatrix@1{{\O_X}\ar[r]^-{d} & {\Omega^1_{X/S}} \ar[r]^-{d} & {\Omega^2_{X/S}}}.$$
	Note that $\Omega^2_{X/S}$ is $R$-torsion, and is supported only on the non-smooth locus of $f$;
	in particular, when $X$ is $S$-smooth we have $\Omega^2_{X/S}=0$.  Simply put, we do not have much to say 
	about the relative de Rham cohomology, as in general it can be somewhat pathological (due to torsion,
	non-degeneration of the Hodge to de Rham spectral sequence, and the presence of $\Omega^2_{X/S}$ which ``morally" should be zero as $X$ has relative dimension one
	over $S$).  For what little we can prove, see Proposition \ref{integraldR} and Remark \ref{lengthomega2}.
\end{remark}

In order to apply certain results of Bloch \cite{Bloch} and Liu-Saito \cite{LiuSaito} in this section, we will need to 
assume that the residue field $k$ of $R$ is perfect, and that $X$ is regular with geometrically connected generic fiber.
Note that although Bloch \cite{Bloch} works only with $R$ equal to the localization at some finite place of the ring of integers in a number field, the proofs
of his results that we use below hold for general $R$ with perfect residue field. 
If the assumption that $k$ is perfect can be removed\footnote{The referee has pointed out to us that the main difficulty with imperfect 
$k$ in \cite{Bloch} and \cite{LiuSaito} is with the definition of the Artin conductor} in \cite{Bloch} and \cite{LiuSaito}, then our results in this section
that depend on these references
(i.e. Lemmas \ref{dRcond}--\ref{dRcondmultfib}, Theorem \ref{dRcondbound}, and Proposition \ref{integraldR})
 will hold without any hypotheses on $k$ (though we do use the assumption that $X_K$ is geometrically connected in the proof of Lemma \ref{dRcond} and elsewhere).
We will also often need to impose the condition
that the greatest common divisor of the multiplicities of the irreducible components of the closed fiber $X_k$ of $X$ is 1;
we will abbreviate this condition by saying that $X_k$ {\em is not a multiple fiber}.  
By perfectness of $k$, this is equivalent to requiring that the greatest common divisor of the {\em geometric} multiplicities
is $1$.  Note that under these hypotheses $X$ is admissible (by Proposition \ref{RaynaudCrit}).

For a prime $\ell$ different from the residue characteristic of $R$, we recall that
the (geometric) {\em Artin conductor of $X$ over $S$} is the positive integer
$$\Art(X/S):=-\chi(X_{\overline{K}})+\chi(X_{\overline{k}})+\Sw H^1_{\et}(X_{\overline{K}}, \Q_{\ell}),$$
where $\chi$ is the usual \'etale Euler characteristic and $\Sw$ is the Swan conductor of the $\ell$-adic
$\Gal(\overline{K}/K)$-representation $H^1_{\et}(X_{\overline{K}}, \Q_{\ell})$.
We can also define the (arithmetic) {\em Artin conductor of $X_K$ over $K$} 
$$\Art(X_K/K):= \dim_K H^1_{\et}(X_{\overline{K}},\Q_{\ell}) - \dim_K H^1_{\et}(X_{\overline{K}},\Q_{\ell})^{I} + \Sw H^1_{\et}(X_{\overline{K}}, \Q_{\ell}),$$
where the vector space in the second term is the subspace of $ H^1_{\et}(X_{\overline{K}},\Q_{\ell})$ fixed by an inertia group $I$ of $\Gal(\overline{K}/K)$
(i.e. \hspace{-1ex}$\Art(X_K/K)$ is the usual Artin conductor of the local $\ell$-adic Galois representation on \'etale cohomology).
Both $\Art(X/S)$ and $\Art(X_K/K)$ are independent of the chosen prime $\ell\neq \Char(k)$, and one has the relation\footnote{Strictly speaking, \cite[Lemma (1.2) (i)]{Bloch}
is proved under the assumption that $k$ is finite and $K$ is a number field.  When $X_k$ is not a multiple fiber, one can use \cite[Proposition 1]{LiuGenus2} instead, which does not impose these hypotheses on $k$ or $K$.  In fact, we will only need (\ref{artinrel}) when $X_k$ is not a multiple fiber.} 
\cite[Lemma (1.2) (i)]{Bloch}
\begin{equation}
	\Art(X_K/K) = \Art(X/S) -(n-1),\label{artinrel}
\end{equation}
where the closed fiber of $X$ has $n$ geometric irreducible components.
Furthermore, if $J_K$ is the Jacobian of $X_K$, then $\Art(X_K/K)$ is the conductor of $J_K$ \cite[\S1.1]{LiuGenus2}.

Consider the map of complexes $\rho:\Omega^{\bullet}_{X/S,\le 1}\rightarrow \omega^{\bullet}_{X/S}$ as in (\ref{rhomap}).  Since
$\rho$ is an isomorphism on $K$-fibers, the hypercohomology of the mapping cone complex $\cone^{\bullet}(\rho)$
is $R$-torsion, so its Euler characteristic
\begin{equation*}
	\chi(\cone^{\bullet}(\rho)) := \sum_{i\in \Z} (-1)^i \len_R \H^i(X,\cone^{\bullet}(\rho))
\end{equation*}
is defined.  A remarkable theorem of Bloch\footnote{Our sign convention for $\Art(X/S)$ is the opposite of Bloch's, but is consistent with that of Liu-Saito \cite{LiuSaito}.
As we noted above, although Bloch assumes that $R$ is the localization at some finite place of the ring of integers in a number field, his proof is valid for general $R$
with perfect residue field. 
} 
\cite[\S2]{Bloch} relates $\chi(\cone^{\bullet}(\rho))$ and $\Art(X/S)$:

\begin{theorem}[Bloch]\label{BlochMain}
	Let $X$ be any regular admissible curve.  Then $\Art(X/S)= -\chi(\cone^{\bullet}(\rho))$.
\end{theorem}

Let us denote by
\begin{equation}
	\xymatrix{
		{\gamma: H^0(X,\Omega^1_{X/S})} \ar[r] & {H^0(X,\omega_{X/S})} 
	}\label{gammap}
\end{equation}
the map on cohomology induced by $\rho$ in degree 1 (i.e.,~the map $H^0(c_{X/S})$, with $c_{X/S}$ as in Proposition \ref{Omegamap}).
In \cite{LiuSaito}, Liu and Saito study the {\em efficient conductor} of $X_K$:
\begin{equation}
	\xymatrix{
		\Effcond(X_K/K) := \len_R \coker (\gamma)
	}.\label{effcondef}
\end{equation}
They show \cite[Lemma 4]{LiuSaito} that $\Effcond(X_K/K)$ is independent of the choice of regular proper model $X$ of $X_K$, 
and (using Theorem \ref{BlochMain}) that when the closed fiber of $X$ is not a multiple fiber one has \cite[Corollary 2]{LiuSaito}
\begin{equation}  
	\Effcond(X_K/K) = \Art(X/S)-\len_R H^1(X,\Omega^1_{X/S})_{\tors}. \label{effcondformula}
\end{equation}
From this, they deduce their main theorem \cite[Theorem 1]{LiuSaito}:

\begin{theorem}[Liu-Saito]\label{LS}
	Let $X$ be a regular proper $S$-curve with $X_k$ not a multiple fiber.  Then
	$$\Effcond(X_K/K)\le \Art(X_K/K).$$
	Furthermore, if $X$ is semistable then this inequality is an equality.
\end{theorem}

We wish to generalize these results to include information about the hypercohomology of the complexes $\Omega^{\bullet}_{X/S,\le 1}$
and $\omega_{X/S}^{\bullet}$.  In order to do this, we first need two lemmas:

\begin{lemma}\label{dRcond}
	  Let $X$ be a regular proper and cohomologically flat $($i.e.~admissible$)$ relative curve over $S$ with geometrically connected
	  generic fiber. Then the natural map 
	\begin{equation*}
		\xymatrix@1{{\H^0(X,\Omega^{\bullet}_{X/S, \le 1})}\ar[r] &  {H^0(X/R)}}
	\end{equation*}
	is an isomorphism, and
	the map $\gamma$ of $(\ref{gammap})$ fits into a natural map of short exact sequences
	\begin{equation}
		\xymatrix@C-8pt{
			0\ar[r] & {H^0(X,\Omega^1_{X/S})} \ar[r]\ar[d]_-{\gamma} & {\H^1(X,\Omega^{\bullet}_{X/S,\le 1})} \ar[r]\ar[d] & 
				{\ker(d: H^1(X,\O_X)\rightarrow H^1(X,\Omega^1_{X/S}))} \ar[r]\ar@{^{(}->}[d] & 0\\
				0\ar[r] & {H^0(X,\omega_{X/S})} \ar[r] & {H^1(X/R)} \ar[r] & {H^1(X,\O_X)}\ar[r] & 0
		}\label{3term}
	\end{equation}
	Furthermore, this diagram depends only on the generic fiber $X_K$ of $X$
	in the following sense: if $X'$ is any other choice of regular proper $($automatically
	cohomologically flat by Remark $\ref{admiexist})$ model of $X_K$, 
	then the diagrams $(\ref{3term})$ corresponding to $X$ and to $X'$ are canonically identified.
\end{lemma}

\begin{proof}
	Since $X$ is proper with geometrically connected generic fiber, we have $f_*\O_X\simeq \O_S$ so the map $d:H^0(X,\O_X)\rightarrow H^0(X,\Omega^1_{X/S})$ is 
	zero.  As $\rho$ is compatible with filtrations, we thus have a commutative diagram
	\begin{equation*}
		\xymatrix{
			{\H^0(X,\Omega^{\bullet}_{X/S,\le 1})} \ar[r]^-{\simeq}\ar[d] & {H^0(X,\O_X)}\ar@{=}[d] \\
			{\H^0(X,\omega_{X/S}^{\bullet})} \ar[r]^-{\simeq} & H^0(X,\O_X)
		}
	\end{equation*}
	giving the first claim. 	Moreover, $\rho$ is also compatible with the usual hypercohomology spectral sequences, so using 
	the first claim we obtain the diagram (\ref{3term}), in which the right vertical map is obviously injective. 
	  To prove the claimed independence of proper regular model, it suffices in view of Theorem \ref{canonicalintegral} and the naturality
	  of the formation of (\ref{3term}) to show that the top row of (\ref{3term})
	depends only on $X_K$.  By the structure theory for regular proper models of $X_K$, we may suppose that
	there is a proper birational map $X'\rightarrow X$ extending the identity on $K$-fibers, and we wish to show that the resulting homomorphism
	\begin{equation}
		\xymatrix@C-10pt{
			0\ar[r] & {H^0(X,\Omega^1_{X/S})} \ar[r]\ar[d] & {\H^1(X,\Omega^{\bullet}_{X/S,\le 1})}\ar[r]\ar[d] & 
			{\ker(H^1(X,\O_{X})\rightarrow H^1(X,\Omega^1_{X/S}))}\ar[r]\ar[d] & 0\\
			0\ar[r] & {H^0(X',\Omega^1_{X'/S})} \ar[r] & {\H^1(X',\Omega^{\bullet}_{X'/S,\le 1})}\ar[r] & 
			{\ker(H^1(X',\O_{X'})\rightarrow H^1(X',\Omega^1_{X'/S}))}\ar[r] & 0
		}\label{indepK}
	\end{equation}
	is an isomorphism.  The left vertical map is an isomorphism by \cite[Lemma 4]{LiuSaito}.
	The right vertical map of (\ref{indepK}) is deduced from the commutative square
	\begin{equation*}
		\xymatrix{
			{H^1(X,\O_{X})} \ar[r]\ar[d] & {H^1(X,\Omega^1_{X/S})}\ar[d] \\
				{H^1(X',\O_{X'})} \ar[r] & {H^1(X',\Omega^1_{X'/S})}
		}
	\end{equation*}
	whose left edge is an isomorphism (as $X$ is regular and $X'\rightarrow X$ is proper birational)
	and whose right edge is injective by \cite[Remark 3]{LiuSaito}.  We conclude that the right vertical map
	of (\ref{indepK}) must be an isomorphism, and hence that the middle vertical map in (\ref{indepK}) is too.
\end{proof}

When the closed fiber of $X$ is not a multiple fiber, then each term in the top row of 
(\ref{3term}) is in fact a submodule of the corresponding term in the bottom row:	
\begin{lemma}\label{dRcondmultfib}	
	  Keep the notation and hypotheses of Lemma $\ref{dRcond}$, and suppose in addition that $X_k$ is not a multiple fiber.
	  Then all three vertical maps in $(\ref{3term})$ are injective, and the canonical map
		\begin{equation*}
			\xymatrix@1{{\H^2(X,\Omega^{\bullet}_{X/S, \le 1})}\ar[r] & H^2(X/R)}
		\end{equation*}
	is surjective with kernel $\H^2(X,\Omega^{\bullet}_{X/S, \le 1})_{\tors}$. 
\end{lemma}

\begin{proof}
	When $X_k$ is not a multiple fiber, 
	 the map $\gamma$ is injective thanks to \cite[Corollary 1]{LiuSaito}, so since
	 the right vertical map of (\ref{3term}) is always injective, the middle vertical map of (\ref{3term}) must be injective as well.
	For the second claim, as in the proof of Lemma \ref{dRcond} we have a commutative diagram
	\begin{equation*}
		\xymatrix{
			{H^1(X,\Omega^1_{X/S})} \ar[r]\ar[d] & {\H^2(X,\Omega^{\bullet}_{X/S,\le 1})} \ar[d] \\
			{H^1(X,\omega_{X/S})} \ar[r]^{\simeq} & {\H^2(X,\omega^{\bullet}_{X/S})} 
		}
	\end{equation*}
	in which the bottom horizontal map is an isomorphism by Proposition \ref{inthodge}.
	The left vertical map of this diagram is surjective by \cite[Corollary 1]{LiuSaito} (again, using the assumption that $X_k$ is not a multiple fiber), 
	so the same is true of the right vertical map.
	We know that the right vertical map is moreover an isomorphism after extending scalars to $K$, and that $\H^2(X,\omega_{X/S}^{\bullet})$
	is torsion free thanks to Proposition \ref{inthodge} and Theorem \ref{grothendieckduality};  our second assertion follows immediately.	  	
\end{proof}

\begin{definition}\label{dRconddef}
	  Let $X/S$ be a regular proper and cohomologically flat relative curve with geometrically connected generic fiber. We
	define the {\em de Rham conductor} of $X_K/K$ to be the integer
	\begin{equation*}
		\dRcond(X_K/K):= \len_R \coker( \H^1(X,\Omega^{\bullet}_{X/S,\le 1})\rightarrow  \H^1(X,\omega^{\bullet}_{X/S})).
	\end{equation*}
	By Lemma \ref{dRcond}, this integer depends only on $X_K$.
\end{definition}

The numerical invariant $\dRcond(X_K/K)$ satisfies the following analogue of (\ref{effcondformula}) and Theorem \ref{LS}:

\begin{theorem}\label{dRcondbound}
	Let $X$ be a regular proper curve over $S$ with geometrically connected generic fiber, and suppose that the closed fiber $X_k$ is not a multiple fiber.  Then
	 \begin{equation}
	 	\dRcond(X_K/K) = \Art(X/S)- \len_R \H^2(X,\Omega^{\bullet}_{X,\le 1})_{\tors}. \label{dRcondformula}
	 \end{equation}
	  Furthermore, we have inequalities
	  \begin{equation}
	  	\Effcond(X_K/K) \le \dRcond(X_K/K) \le \Art(X_K/K),\label{ineqs}
	  \end{equation}
	  and when $X$ is semistable, each of these inequalities is an equality.
\end{theorem}

\begin{proof}
	First note that the final assertion concerning the case of semistable $X$ follows immediately from Theorem \ref{LS},
	given the inequalities (\ref{ineqs}).  Thus, we need only establish (\ref{dRcondformula}) and (\ref{ineqs}).
	
	  As the hypotheses and conclusions are unaffected by base change to the strict henselization of $R$, we may
	  assume that $k$ is separably closed, and hence algebraically closed as $k$ is perfect.
	For convenience, we will abbreviate $\Omega^{\bullet}_{X/S,\le1}$ and $\omega_{X/S}^{\bullet}$ 
	simply by $\Omega^{\bullet}$ and $\omega^{\bullet}$, respectively.
	The map $\rho:\Omega^{\bullet} \rightarrow \omega^{\bullet}$ fits into an exact triangle
	\begin{equation*}
		\xymatrix{
  			 {\Omega^{\bullet}} \ar[r]^-{\rho} & {\omega^{\bullet}}\ar[r] & {\cone^{\bullet}(\rho)}
		},
	\end{equation*}
	from which we deduce the usual long exact sequence in hypercohomology.  Using Lemmas \ref{dRcond} and \ref{dRcondmultfib},
	we extract from this the exact sequence
	\begin{equation}
			\xymatrix{
			0 \ar[r] & {\H^1(\Omega^{\bullet})} \ar[r] & {\H^1(\omega^{\bullet})} \ar[r] & {\H^1(\cone^{\bullet}(\rho))} \ar[r] &
			  {\H^2(\Omega^{\bullet})} \ar[r]^-{\H^2(\rho)} & {\H^2(\omega^{\bullet})} \ar[r] & 0
			}\label{coneseq}
	\end{equation}
	(with self-evident abbreviated notation), and that one has $\H^i(\cone^{\bullet}(\rho))=0$ for $i\neq 1$.
	In particular, we have $\chi(\cone^{\bullet}(\rho)) = -\len_R \H^1(\cone^{\bullet}(\rho))$, so the
	exact sequence (\ref{coneseq})	tells us that
	\begin{equation*}
		\dRcond(X_K/K) = -\chi(\cone^{\bullet}(\rho)) -\len_R \ker(\H^2(\rho))
	\end{equation*}
	Applying Theorem \ref{BlochMain} and using Lemma \ref{dRcondmultfib} again, we deduce (\ref{dRcondformula}).
	
	Now the first inequality in (\ref{ineqs}) easily follows from the kernel-cokernel exact sequence arising from the diagram (\ref{3term}) 
	and the definition of the efficient and de Rham conductors.
	In order to establish the second inequality, we proceed as follows.
	Suppose that the irreducible components of the closed fiber $X_k$ are $C_1$, $C_2,$\ldots, $C_n$.
	Thanks to (\ref{artinrel}), to prove the second inequality it suffices to show that $\len_R \H^2(\Omega^{\bullet})_{\tors} \ge n-1$.  We do this by refining
	the argument in the proof of \cite[Theorem 1]{LiuSaito}.
	Observe that 
	\begin{equation*}
		\len_R \H^2(X,\Omega^{\bullet})_{\tors} \ge \dim_k \H^2(X,\Omega^{\bullet})_{\tors}\otimes_R k =
		\dim_k (\H^2(X,\Omega^{\bullet})\otimes_R k) -1,
	\end{equation*}
	where the equality of the second and third terms results from the fact that  the $K$-vector space $\H^2(X,\Omega^{\bullet}_{X/S,\le 1})\otimes K$
	is one-dimensional.  
	It is not hard to see that the functor $\H^2(X,\cdot)$ is right exact on two-term complexes of coherent sheaves, so the
	the natural base change map 
	\begin{equation}
		\xymatrix{
			{\H^2(X,\Omega^{\bullet})\otimes_R k} \ar[r] & {\H^2(X_k,\Omega^{\bullet}_{X_k})}
	}\label{bc}
	\end{equation}
	is surjective.\footnote{It is even an isomorphism.}  Thus, we will be done if we can show that 
	$\dim_k \H^2(X_k,\Omega^{\bullet}_{X_k}) \ge n$.
	Give $C_i$ its reduced structure.
	  Associated to the natural maps $\iota_i: C_i\rightarrow X_k$ is
	the (filtered) map of two-term filtered complexes 
	\begin{equation}
		\Omega^{\bullet}_{X_k}\rightarrow \bigoplus_i {\iota_i}_* \Omega^{\bullet}_{C_i/k},\label{componentstrick}
	  \end{equation}
	  which induces a commutative diagram
	\begin{equation*}
		\xymatrix{
			{H^1(X_k,\Omega^1_{X_k/k})}\ar@{->>}[r]\ar[d] & {\H^2(X_k,\Omega^{\bullet}_{X_k})}\ar[d] \\
			{\bigoplus_i H^1({C}_i, \Omega^1_{{C}_i/k})}\ar@{->>}[r] &
			{\bigoplus_i \H^2({C}_i, \Omega^{\bullet}_{{C}_i/k})}
		}
	\end{equation*}
	 in which the horizontal maps are surjective (due to the vanishing of $H^2$ of coherent sheaves on curves over a field).
	Now the left vertical map is also surjective,
	since the cokernel of (\ref{componentstrick}) is supported on finitely many closed points
	and $H^1$ is right exact and annihilates skyscraper sheaves.  It follows that the right vertical map is surjective as well.  
	Provided that each $\H^2({C}_i, \Omega^{\bullet}_{{C}_i/k}) = \H^2_{\dR}({C}_i/k)$ is nonzero,
	we could conclude that $\dim_k \H^2(X_k,\Omega^{\bullet}_{X_k})\ge n$, so in view of the surjectivity of (\ref{bc}) we would be done.
	
	It remains to show that $H^2_{\dR}(C'/k)\neq 0$ for any proper $k$-scheme $C'$ of dimension 1.  
	  We can adapt the argument using (\ref{componentstrick}) to show that $H^2_{\dR}(C'/k)$ maps onto $\bigoplus_j H^2_{\dR}(C'_j/k)$,
	where $\{C'_j\}$ are the normalizations of the irreducible components of $C'$.  But the Hodge to de Rham spectral sequence
	for each smooth curve $C_j'$ degenerates, so $H^2_{\dR}(C'_j/k)$ is 1-dimensional over $k$ for each $j$.
\end{proof}

\begin{remark}
	In fact, we do not need the assumption in Definition \ref{dRconddef} that $X/S$ be cohomologically flat.  
	Certainly, the map $\rho$ of (\ref{rhomap}) always yields  map 
	 $$\H^1(X,\Omega^{\bullet}_{X/S,\le 1})\rightarrow  \H^1(X,\omega^{\bullet}_{X/S})$$ on hypercohomology,
	 and we can always contemplate the length of the cokernel of this map (even in the absence of the regularity hypothesis).
	 The point is that an analogue of Lemma \ref{dRcond} holds when $X$ is only assumed to be regular, and in this way
	 one shows that $\dRcond(X_K/K)$ really only depends on $X_K$.
\end{remark}

\begin{remark}
	Let $X$ be a regular proper curve over $S$ with geometrically connected generic fiber and suppose that the closed fiber $X_k$ is not a multiple fiber. 
	It is natural to ask for a conceptual interpretation of the de Rham conductor, relative to the efficient and Artin conductors.
	In fact, the nonnegative numerical invariants $\dRcond(X_K/K)-\Effcond(X_K/K)$ and $\Art(X_K/K)-\dRcond(X_K/K)$
	give a quantitative measure of the degree to which the Hodge to de Rham spectral sequence for the truncated de Rham
	complex $\Omega_{X/S,\le 1}^{\bullet}$ degenerates.  More precisely, using the construction of (\ref{3term}) and Lemma \ref{dRcondmultfib},
	one has a short exact sequence of $R$-modules
	\begin{equation}
		\xymatrix{
			0\ar[r] & {\Im(d)} \ar[r] & {H^1(X,\Omega^1_{X/S})_{\tors}} \ar[r] & {\H^2(X,\Omega^{\bullet}_{X/S,\le 1})_{\tors}}\ar[r] & 0
		}\label{imd}
	\end{equation}
	where $d$ is the edge map 
	\begin{equation*}
		\xymatrix{
			{d:H^1(X,\O_X)}\ar[r] & {H^1(X,\Omega^1_{X/S})} 
			}
	\end{equation*}
	induced by the usual derivation (which must have 
	torsion image due to the degeneration of the Hodge to de Rham spectral sequence of $X_K$ over $K$).  
	Subtracting (\ref{effcondformula}) from (\ref{dRcondformula}) and using (\ref{imd}), we deduce the equality
	\begin{equation}
		\dRcond(X_K/K) - \Effcond(X_K/K) = \len_R(\Im(d)).\label{dediff}
	\end{equation}
	In particular, the de Rham and efficient conductors of $X_K$ over $K$ are equal if and only if the Hodge to de Rham 
	spectral sequence for the truncated de Rham complex $\Omega^{\bullet}_{X/S,\le 1}$ degenerates.
	In a similar manner, if $n$ denotes the number of geometric irreducible components of $X_k$,
	the proof of Theorem \ref{dRcondbound} and (\ref{imd}) show that 
	\begin{equation}
		\Art(X_K/K) - \dRcond(X_K/K) = \len_R(H^1(X,\Omega)_{\tors}/\Im(d)) - (n-1),\label{addiff}
	\end{equation}
	and this quantity is nonnegative due to (\ref{ineqs}).  Thus, the de Rham and Artin conductors are equal if 
	and only if $\Im(d)$ is as large as theoretically possible.  Note that the right sides of (\ref{dediff})
	and (\ref{addiff}) are invariants of the generic fiber $X_K$, because the left sides are.
\end{remark}

We end this section with a proposition that elucidates the structure of the relative de~Rham cohomology of $X$ over $S$.
As we have indicated, integral de Rham cohomology is often somewhat pathological, so it should come as no surprise that there
is little we can say about it in full generality.  Nonetheless: 

\begin{proposition}\label{integraldR}
	Let $X$ be a proper relative curve over $S$ with geometrically connected generic fiber. Then there are canonical isomorphisms
	\begin{equation*}
		H^0_{\dR}(X/S) \simeq H^0(X,\O_X)\quad\text{and}\quad H^3_{\dR}(X/S) \simeq \frac{\D H^1(X,\Omega^2_{X/S})}{\D d H^1(X,\Omega^1_{X/S})},
	\end{equation*}
	as well as an injective map
	\begin{equation*}
		\xymatrix{
			{H^1_{\dR}(X/S)} \ar@{^{(}->}[r]^-{\alpha} & {\H^1(X,\Omega^{\bullet}_{X/S,\le 1})}
		}
	\end{equation*}
	with 
	$$\len_R \coker (\alpha) \le \len_R H^0(X,\Omega^2_{X/S}).$$
	Moreover, the cokernel of the canonical map $H^0(X,\Omega^2_{X/S})\rightarrow H^2_{\dR}(X/S)$ is naturally isomorphic to
	\begin{equation*}
		\frac{\D \ker (H^1(X,\Omega^1_{X/S}) \xrightarrow{d} H^1(X,\Omega^2_{X/S}) )}{\D dH^1(X,\O_{X})}.
	\end{equation*}
	In particular, if $X$ is regular and $X_k$ is not a multiple fiber then $H^1_{\dR}(X/S)$ is torsion-free,  and if $X_k$ is in addition generically smooth,
	then $H^3_{\dR}(X/S)=0$. 
\end{proposition}

\begin{proof}
	A straightforward analysis of the Hodge to de Rham 
	spectral sequence yields the claimed identifications of $H^0_{\dR}$ and $H^3_{\dR}$, as well
	as the statements concerning $H^2_{\dR}$ and $\alpha$.
	  When $X$ is regular and $X_k$ is not a multiple fiber, we may apply  Lemma \ref{dRcondmultfib} to deduce that we
	have an injection $H^1_{\dR}(X/S)\hookrightarrow H^1(X/R)$; the latter module is torsion free by Proposition \ref{inthodge}.
	When $X_k$ is generically smooth, $\Omega^2_{X/S}$ is a skyscraper sheaf (as in this case the non-smooth locus of $X$ over $S$ consists
	of a finite number of closed points in the closed fiber), so $H^1(X,\Omega^2_{X/S})=0$.
\end{proof}

\begin{remark}\label{lengthomega2}
	It is natural to ask if one can give a bound for $L:=\len H^0(X,\Omega^2_{X/S})$.
	Using the results of this section, it is not hard to show that when $X$ is regular and $X_k$ is not a multiple fiber,
	we have $L \ge n-1$, where $X_k$ has $n$ geometric irreducible components.  When $X_k$ is in addition generically smooth,
	one can also show that $L \le \Art(X/S)$.   
\end{remark}

\section{Base change}\label{BaseChng}
	  Let $X_K$ be a smooth proper and geometrically connected curve over $K$ with an admissible model $X$ over $R$,
	and fix a local extension of discrete valuation rings $R\rightarrow R'$ with associated fraction field extension $K\rightarrow K'$.
	Extension of scalars gives an integral structure $H^1(X/R)\otimes_R R'$ on $H^1_{\dR}(X_{K'}/K')$.  When $X_{K'}$
	has an admissible model $X'$ over $R'$, it is natural to ask how this lattice relates to the canonical integral structure
	$H^1(X'/R')$.	
	  If the residue field extension $k\rightarrow k'$ is separable and the relative ramification index $e(R'/R)$ is 1,
	then $X'=X\times_R {R'}$ is an admissible model of $X_{K'}$, as we noted in the discussion preceding Proposition \ref{BC},
	and hence $H^1(X/R)\otimes_R R'$ and $H^1(X'/R')$ are canonically isomorphic.  This includes the cases $R'=\widehat{R}$
	and $R'=R^{\mathrm{sh}}$, so the most interesting cases have $R$ complete and $[K':K]$ finite.
	For general base change, these two lattices need not coincide, nor must either lattice necessarily contain the other.  
	Nonetheless, we can estimate how far apart these lattices inside $H^1_{\dR}(X_{K'}/K')$ are in general:
	
	\begin{theorem}	\label{latticeposition}
		Let $X_K$ be a smooth, proper and geometrically connected curve over $K$ with an admissible model $X$ over $R$,
		and assume that $[K':K]$ is finite and that $X_{K'}$ has an admissible model $X'$ over $R'$.  Set  
		$$L:=H^1(X/R)\otimes_R R'\quad\text{and}\quad L':=H^1(X'/R'),$$
		considered as $R'$-lattices inside $H^1_{\dR}(X_{K'}/K')$.   Denote by $e$ the relative ramification index of $R'$ over $R$,
		 fix a uniformizer $\pi$ of $R$, and put $E:=\Effcond(X_K/K)$ $($see $(\ref{effcondef})$ for the definition$)$.
		Then
		\begin{enumerate}
			\item $\pi^EL\subseteq L'$ and $\pi^E L'\subseteq L$.\label{latpos1}
			\item Each of $L/L\cap L'$ and $L'/L\cap L'$ is of finite $R'$-length at most $e E$.\label{latpos2}
		\end{enumerate}
	\end{theorem}

\begin{proof}
	As usual, put $S=\Spec R$ and $S'=\Spec R'$.
	Thanks to Theorem \ref{canonicalintegral} and Proposition \ref{BC}, we may suppose that $X$ is regular and that $X'$ is the normalization
	of $X_{S'}:=X\times_S S'$.  The normalization map $\rho:X'\rightarrow X_{S'}$ is then finite and birational, so we have natural maps
	$\O_{X_{S'}}\rightarrow \rho_* \O_{X'}$ and $\rho_*\omega_{X'/S'}\rightarrow \omega_{X_{S'}/S'}$ 
	whose cokernels are supported on the closed fiber of $X_{S'}$.  These maps fit naturally into diagrams of 2-term (horizontal) complexes on $X_{S'}$

	\begin{minipage}{3in}
	\begin{equation}
		\xymatrix{
			{\O_{X_{S'}}} \ar[r]^-{d}\ar[d] & {\omega_{X_{S'}/S'}} \ar@{=}[d] \\
			{\rho_*\O_{X'}} \ar[r]\ar@{=}[d] & {\omega_{X_{S'}/S'}} \\
			{\rho_*\O_{X'}} \ar[r]_-{d} & {\rho_*\omega_{X'/S'}}\ar[u]
		}\label{D1}
	\end{equation}
	\end{minipage}
	\hfill
	\begin{minipage}{3in}
	\begin{equation}
		\xymatrix{
			{\O_{X_{S'}}} \ar[r]^-{d}\ar@{=}[d] & {\omega_{X_{S'}/S'}}  \\
			{\O_{X_{S'}}} \ar[r]\ar[d] & {\rho_*\omega_{X'/S'}}\ar[u] \\
			{\rho_*\O_{X'}} \ar[r]_-{d} & {\rho_*\omega_{X'/S'}}\ar@{=}[u]
		}\label{D2}
	\end{equation}
	\end{minipage}
	where in each case the middle horizontal arrow is deduced from the bottom one, and the diagrams commute as this may be checked
	over the generic fiber, where it is obvious.  We claim that for each 2-term (horizontal) complex occurring in (\ref{D1}) and (\ref{D2}),
	the usual ``Hodge to de Rham" spectral sequence for hypercohomology degenerates at the $E_1$-stage.  Due to such degeneration
	on the generic fiber, as in the proof of  Proposition \ref{inthodge} it suffices to show that $H^i(X_{S'}, \omega_{X_{S'}/S'})$
	and $H^i(X_{S'}, \rho_*\omega_{X'/S'})$ are torsion free for $i=0,1$.  This follows easily from Theorem \ref{grothendieckduality}
	(using that $X$ and $X'$ are admissible) after making the canonical identifications
	\begin{eqnarray*}
	 H^i(X_{S'}, \omega_{X_{S'}/S'})\simeq H^i(X,\omega_{X/S})\otimes_R R' & \text{and} &
	H^i(X_{S'}, \rho_*\omega_{X'/S'})\simeq H^i(X',\omega_{X'/S'})
	\end{eqnarray*}
	that result from the compatibility of the dualizing sheaf with base change and the finiteness of $\rho$,
	respectively.  We deduce from (\ref{D1}) and (\ref{D2}) the following commutative diagrams of short exact sequences of free $R'$-modules:
	\begin{equation}
		\xymatrix{
			0\ar[r] & {H^0(X,\omega_{X/S})\otimes_R R'} \ar[r]\ar@{=}[d] & {H^1(X/R)\otimes_R R'} \ar[r]\ar[d]^-{i} & {H^1(X,\O_X)\otimes_R R'} \ar[r]\ar[d]^-{\alpha} & 0\\
			0\ar[r] & {H^0(X,\omega_{X/S})\otimes_R R'} \ar[r] & {\H^1(\rho_*\O_{X'}\rightarrow \omega_{X_{S'}/S'})} \ar[r] & {H^1(X',\O_{X'})} \ar[r] & 0\\
			0\ar[r] & {H^0(X',\omega_{X'/S'})} \ar[r]\ar[u]^-{\beta} & {H^1(X'/R')} \ar[r]\ar[u]^-{i'} & {H^1(X',\O_{X'})} \ar[r]\ar@{=}[u] & 0
		}\label{C1}
	\end{equation}
	\begin{equation}
		\xymatrix{
			0\ar[r] & {H^0(X,\omega_{X/S})\otimes_R R'} \ar[r] & {H^1(X/R)\otimes_R R'} \ar[r] & {H^1(X,\O_X)\otimes_R R'} \ar[r]\ar@{=}[d] & 0\\
			0\ar[r] & {H^0(X',\omega_{X'/S})} \ar[r]\ar[u]^-{\beta} & {\H^1(\O_{X_{S'}}\rightarrow \rho_*\omega_{X'/S'})}\ar[u]^-{j}\ar[d]^-{j'} \ar[r] & 
			{H^1(X,\O_{X})\otimes_R R'} \ar[r]\ar[d]^-{\alpha} & 0\\
			0\ar[r] & {H^0(X',\omega_{X'/S'})} \ar[r]\ar@{=}[u] & {H^1(X'/R')} \ar[r] & {H^1(X',\O_{X'})} \ar[r] & 0
		}\label{C2}
	\end{equation}
	Since $\alpha$ and $\beta$ have torsion-free sources and become isomorphisms after extending scalars to $K'$,
	we conclude that both maps are injective and have finite-length cokernels.   
	  The diagrams (\ref{C1}) and (\ref{C2}) then imply that $i$ and $j$ are injective, and that (via the natural maps)
	\begin{eqnarray}
		\coker i \simeq \coker \alpha \simeq \coker j' &\text{and}& \coker i'\simeq \coker \beta \simeq \coker j.\label{cokerisoms}
	\end{eqnarray}
	  On the other hand, Grothendieck duality for $\rho$ ensures that the maps $\alpha$ and $\beta$ are dual, and that there is a 
	  canonical isomorphism of finite length $R'$-modules
	\begin{equation*}
		\coker \beta \simeq \Ext^1_{R'}(\coker \alpha,R').
	\end{equation*} 
	In particular $\coker \alpha$ and $\coker \beta$ have the same $R'$-length. 
	Thus, Theorem \ref{latticeposition} follows easily from (\ref{cokerisoms}), the diagrams (\ref{C1}) and (\ref{C2}), and the bound
	\begin{equation}
		\len_{R'}(\coker \beta) \le eE.\label{lastbound}
	\end{equation}
	To prove (\ref{lastbound}), we first claim that $\beta$ fits into a commutative diagram
	\begin{equation}
	\xymatrix{
		{H^0(X,\Omega^1_{X/S})\otimes_R R'} \ar[r]\ar[d]_-{\gamma_{R'}} & {H^0(X',\Omega^1_{X'/S'})\ar[d]} \\
		{H^0(X,\omega_{X/S})\otimes_R R'} & \ar[l]^-{\beta} {H^0(X',\omega_{X'})}
		}\label{cooltrick}
	\end{equation}
	where the vertical maps are induced by the canonical maps from differentials to dualizing sheaves 
	  (i.e. the map $\gamma_{R'}$ is the base change of (\ref{gammap}); {\em cf.} also Proposition \ref{Omegamap})
	and the top horizontal map is the usual pullback map on differentials associated to $\rho: X'\rightarrow X_{S'}$.
	Indeed, since $H^0(X,\omega_{X/S})\otimes_R R'$ is torsion free, the commutativity of (\ref{cooltrick})
	may be checked after extending scalars to $K'$, where it follows immediately from the construction of the maps involved.
	As a particular consequence of the commutativity of (\ref{cooltrick}), we deduce that 
	\begin{equation*}
	\len_{R'}( \coker \beta) \le \len_{R'} (\coker \gamma_{R'}).
	\end{equation*}
	On the other hand, by the definition (\ref{effcondef}) of $E:=\Effcond(X_K/K)$ and the fact that $\gamma_{R'}$ is obtained from the map $\gamma$
	of (\ref{gammap}) by base change, we have
	\begin{equation*}
		\len_{R'}( \coker \gamma_{R'})  = e\len_{R}( \coker \gamma) = eE,
	\end{equation*}
	and the required bound (\ref{lastbound}) follows at once.
	\end{proof}

\begin{remarks}\label{latcontainment}
			As we observed immediately preceding Proposition \ref{BC}, the hypotheses of Theorem \ref{latticeposition} are satisfied by any smooth and geometrically 			connected curve $X_{K}$ equipped with a rational point over a tamely ramified extension of $K$.
			
			 Note that the bounds given in Theorem \ref{latticeposition} (\ref{latpos1}) do not depend on $R'$ at all.
			We do not know if these bounds are sharp; however, it is easy to see that if one of $L$, $L'$ contains the other,
			then they must be equal (as lattices inside $H^1_{\dR}(X_{K'}/K')$).  Indeed, suppose, for example, that $L\subseteq L'$.
			Then from either of the diagrams (\ref{C1}) or (\ref{C2}), we deduce a commutative diagram
		\begin{equation*}
		\xymatrix{
			0\ar[r] & {H^0(X,\omega_{X/S})\otimes_R R'} \ar[r] & {H^1(X/R)\otimes_R R'}\ar@{^{(}->}[d]\ar[r] & {H^1(X,\O_X)\otimes_R R'} \ar[r]\ar@{^{(}->}[d] & 0\\
			0\ar[r] & {H^0(X',\omega_{X'/S'})} \ar[r]\ar@{^{(}->}[u] & {H^1(X'/R')} \ar[r] & {H^1(X',\O_{X'})} \ar[r] & 0
		}
	\end{equation*}
			from which it follows (by a simple diagram chase) that the left vertical injection is an isomorphism.  By Grothendieck duality for $\rho:X'\rightarrow X_{R'}$, 
			the right map is also an isomorphism, so the same is true for the middle map.  The argument in the case of the reverse inclusion of lattices is similar.
			
			Supposing that $k$ is perfect,  when the closed fiber of $X$ is not a multiple fiber (for example, when $X_{K}$ has a rational point over an unramified extension
			of $K$), we may apply Theorem \ref{LS} of Liu-Saito to bound $E$ above by $\Art(X_K/K)$.
\end{remarks}

\bibliography{mybib}
\end{document}